
\magnification = 1200

\def\authors#1{\halign{\quad##\hfil&&\quad\qquad\qquad\qquad##\hfil\cr#1\crcr}}

\raggedbottom  
\input amssym.def
\input amssym.tex

\font\bigbold = cmb10 scaled \magstep2
\font\bigtenrm = cmr10 scaled \magstep1

\def\bwedge{\bigwedge}

\def\Dom{\rm{Dom}}
\def\Max{\rm{Max}}

\def\Min{\rm{Min}}

\def\pr{\rm{pr}}

\def\club{\clubsuit}
\def\cf{\rm{cf}}

\def\rest{\upharpoonright}  
\def\eop{$\bigstar$}  
\def\1st{1$^{\hbox{st}}$}
\def\2nd{2$^{\hbox{nd}}$}
\def\card#1{\vert{#1}\vert} 

\def\cf{\mathop{\rm cf}}
\def\ref{\mathop{\rm ref}}  
\def\nacc{\mathop{\rm nacc}}
\def\acc{\mathop{\rm acc}}
\def\implies{\Longrightarrow}



\def\deq{\buildrel{\rm def}\over =}
\def\otp{\mathop{\rm otp}}
\def\Min{\mathop{\rm Min}}
\def\cl{\mathop{\rm cl}}
\def\to{\longrightarrow}
\def\mod{\mathop{\rm mod}}

\def\DD{{\cal D}}
\def\FF{{\cal F}}
\def\II{{\cal I}}
\def\PP{{\cal P}}



\def\today{\ifcase\month\or
   January\or February \or March \or April \or May \or June \or
   July \or August \or September \or October \or November \or December \fi
   \space \number\day, \number\year}


\centerline{{\bigbold
Saturated filters at successors of singulars,
weak reflection and}}

\centerline{{\bigbold yet another weak club principle}}


\bigskip
\bigskip
\bigskip

\def\pomak{\negthinspace\negthinspace\negthinspace\negthinspace\negthinspace
\negthinspace\negthinspace\negthinspace\negthinspace\negthinspace}

\centerline{\pomak by \footnote{${}^1$}{\sevenrm Authors partially
supported by the Basic Research Foundation Grant number 0327398
administered by the Israel Academy of Sciences and Humanities.
For easier future reference, note that
this is publication [DjSh 545] in Shelah's bibliography.
The research for this paper was done in the period between January and
August 1994.
The most sincere thanks to James Cummings for his numerous comments
and corrections,
and to Uri Abraham, Menachem Magidor and other 
participants of the Logic Seminar in Jerusalem, Fall 1994,
for their interest. The first author also wishes to thank the Hebrew
University and the Lady Davis Foundation for the Forchheimer Postdoctoral
Fellowship.
}}

\medskip

\authors{{\bigtenrm Mirna D\v zamonja}&{\bigtenrm Saharon
Shelah} 
\cr
University of Wisconsin-Madison&   
 Hebrew University of Jerusalem\cr
Department of Mathematics&Institute of Mathematics\cr
Madison, WI 53706&91904 Givat Ram\cr
USA&Israel\cr
{\sevenrm dzamonja@math.wisc.edu}&{\sevenrm shelah@math.huji.ac.il}\cr}

\vskip 2 true cm

\centerline{\bf Abstract}

\smallskip

Suppose that $\lambda$ is the successor of a singular cardinal
$\mu$ whose cofinality is an uncountable cardinal $\kappa$.
We give a sufficient condition that the club filter
of $\lambda$ concentrating on the points of cofinality $\kappa$ is
not $\lambda^+$-saturated.\footnote{${}^2$}{\sevenrm added in proof:
M. Gitik and S. Shelah have subsequently and by a
different technique shown that the club
filter on such $\lambda$ is never saturated.} The condition is phrased in terms
of
a notion that we call weak reflection. We discuss various properties
of weak reflection.

We introduce a weak version of the $\club$-principle, which we call
$\club^\ast_-$, and show that if it holds on a stationary subset $S$
of $\lambda$, then no normal filter on $S$ is $\lambda^+$-saturated.
Under the above assumptions, $\club^\ast_-(S)$ is true for any
stationary subset $S$ of $\lambda$ which does not contain points
of cofinality $\kappa$. For stationary sets $S$ which concentrate
on points of cofinality $\kappa$, we show that $\club^\ast_-(S)$
holds modulo an ideal obtained through the weak reflection.

\vfill
\eject
\baselineskip=24pt

{\bf \S 0. Introduction.}
Suppose that $\lambda=\mu^+$ and $\mu$ is an infinite cardinal of
cofinality $\kappa$.
We revisit the classical question of whether a normal
filter on $\lambda$ can be $\lambda^+$-saturated.
We are in particular concerned with the case of $\kappa$
uncountable and less than $\mu$. We are mostly interested in the club filter
on $\lambda$.

While the richness of the literature on the subject provides us with
a strong motivation for a further study, it also prevents
us from giving a complete history and bibliography
involved. We shall give a list
of those references which are most directly connected
or used in our results,
and for further reading we can suggest looking at the 
references mentioned in the papers that we refer to.

It is well known that for no regular $\theta>\aleph_0$ the club 
filter on $\theta$ can be $\theta$-saturated, but
modulo the existence of huge or other large cardinals it is consistent
that $\aleph_1$ carries an $\aleph_2$-saturated normal filter (Kunen [Ku1]),
or in fact that any uncountable
regular $\aleph_\alpha$ carries an $\aleph_{\alpha+1}$-saturated
normal filter (Foreman [Fo]). In these arguments
the 
saturated filter obtained is not the club filter. 
As a particular case of [Sh 212, 14] or [Sh 247, 6],
if $\sigma=\rho^+$, then
the club filter $\DD_\sigma$ restricted to the elements
of $\sigma$ of a fixed cofinality $\theta\neq\cf(\rho)$
is not $\sigma^+$-saturated.
It is consistent that $\DD_{\aleph_1}$ is $\aleph_2$-saturated, as is shown in
[FMS] and also in [SvW].

If $\sigma<\rho$ and $\rho$ is a regular cardinal, we use $S_\rho^\sigma$ to 
denote the set of elements of $\sigma$ which have cofinality $\rho$.
Let $\lambda$ be as above.

In the first section of the present paper, we give a sufficient condition
that $\DD_\lambda\rest S^\lambda_\kappa$ is not $\lambda^+$-saturated.
Here $\aleph_0<\kappa=\cf(\mu)<\mu$.
The condition is a reflection property, which we shall call
{\it weak reflection}, as we show that it is weaker than some known
reflection properties. We discuss the properties of weak reflection in more
detail in \S1. Of course, the main part of the section is to show that the
appropriate form of the weak reflection
indeed suffices for $\DD_\lambda\rest S^\lambda_\kappa$ not
to be $\lambda^+$-saturated, which is done in 1.13. In 1.15 the argument is
generalized to some other normal filters on $\lambda$, and [Sh 186,\S3]
is revisited.

In the second section
we introduce the combinatorial principle $\club^\ast_-$
which 
has the property that, if $\club^\ast_-(S)$ holds, then no normal filter
on $S$ is $\lambda^+$-saturated. Here $S$ is a stationary subset of $\lambda$.
The $\club^\ast_-$ is a weak form of $\club$. The $\club$-principle
was first introduced for $\aleph_1$ in Ostaszewski [Os], and later investigated
in a more general setting in [Sh 98] and elsewhere.

If $\lambda$ is the successor of the singular cardinal $\mu$ whose
cofinality is $\kappa$, then $\club^\ast_-(\lambda\setminus S^\lambda_\kappa)$
is true just in $ZFC$.
As a corollary
of this we obtain
an alternative proof of a part of the result of [Sh 212, 14] or [Sh 247, 6]. 
On the other hand, $\club^\ast_-(S^\lambda_\kappa)$ only holds
modulo an ideal defined through the weak reflection (in \S1), so
this gives a connection with the results of the first section.
For the case of $\mu$ being a strong limit, $\club^\ast_-$ just
becomes the already known $\club^\ast$.
Trying to apply
to $\club^\ast_-$ arguments which work for the corresponding version
of $\diamondsuit$,
we came up with a question we could not answer,
so we pose it in 2.6.

Before we proceed to present our results, we shall introduce some notation
and conventions that will be used throughout the paper. 

\medskip
 
{\bf Notation 0.0}(0)
Suppose that $\gamma\ge\theta$ and $\theta$ is a regular
cardinal.
Then 
$$S_{<\theta}^\gamma=\{\delta <\gamma:\aleph_0\leq\cf(\delta) <\theta\},
\quad\hbox{ if }\theta>\aleph_0.$$
$$S_{\theta}^\gamma=\{\delta <\gamma: \cf(\delta)=\theta\}.$$

More generally, we use
$S^\gamma_{{\rm r}\, \theta}$ for ${\rm r}\in
\{<,\le,=,\neq,>,\ge\}$ to describe
$$S^\gamma_{{\rm r}\,\theta}=\{\delta<\gamma:\,\aleph_0\le\cf(\delta)\,\&\,
\cf(\delta)\,{\rm r}\,\theta\}.$$

(1) $SING$ denotes the class of
singular ordinals, that is, all ordinals $\delta$ with $\cf(\delta)<\delta$.

$LIM$ denotes the class of all limit ordinals.

(2) For $\lambda$ a cardinal with $\cf(\lambda)>\aleph_0$,
we denote by $\DD_\lambda$
the club filter on $\lambda$.
The ideal of non-stationary subsets of $\lambda$ is denoted by
$J_\lambda$.

(3) If $C\subseteq\lambda$, then
$$\acc(C)=\{\alpha\in C:\,\alpha=\sup(C\cap\alpha)\}\quad\hbox{and}\quad
\nacc(C)=C\setminus\acc(C).$$

\medskip

We now go on to the first section of the paper.

\bigskip

{\bf \S1. Saturated filters on
the successor of a singular cardinal
of uncountable cofinality.}
Suppose that $\mu$ is a cardinal with the property
$$\mu>\cf(\mu)=\kappa>\aleph_0,$$
and that $\lambda=\mu^+$. We wish to discuss the 
saturation of $\DD_\lambda\rest S^\lambda_\kappa$
and some other
normal filters on $\lambda$. For the reader's
convenience, we recall some relevant definitions and notational
conventions.
The cardinals $\lambda,\mu$ and $\kappa$ as above will be fixed throughout
this section.

\smallskip

{\bf Notation 1.0}(0) Suppose that $\DD$ is a filter on the set $A$. The
{\it dual ideal\/} $\II$ of $\DD$ is defined as
$$\II=\{a\subseteq A:\,A\setminus a\in \DD\}.$$
A set $a\subseteq A$ is $\DD$-{\it positive\/} 
or
$\DD$-{\it stationary\/}, if $a\notin \II$.
The family of all $\DD$-stationary sets
is denoted by ${\DD}^+$. A subset of $A$ which is not $\DD$-stationary,
is referred to as $\DD$-{\it non-stationary\/}.

(1) Suppose that $\sigma$ is a cardinal. 
A filter $\DD$ on a set
$A$ is $\sigma$-{\it saturated\/} iff  there are no $\sigma$
sets which are all
$\DD$-stationary, but no
two of them have a $\DD$-stationary intersection.
If we denote the dual ideal of $\DD$ by $\II$, then 
this is equivalent to saying that
$\PP(A)/\II$ has $\sigma$-ccc.

For a cardinal $\rho$, the filter $\DD$ is $\rho$-{\it complete\/}
if it is closed under taking intersections of $<\rho$ of its elements.

(2) A $\sigma$-complete
filter $\DD$ on a cardinal $\sigma$ is {\it normal\/}, if
for any sequence $X_\alpha(\alpha<\sigma)$ of elements of $\DD$,
the diagonal intersection
$$\Delta_{\alpha<\sigma}X_\alpha\deq\{\beta\in\sigma:\,\alpha<\beta
\implies\beta\in X_\alpha\},$$
is an element of $\DD$, and $\DD$ contains all final segments of $\sigma$.

(3) Suppose that $\DD$ is a filter on $A$ and $S\subseteq A$ is $\DD$-stationary.
We use $\DD\rest S$ to denote
$$\DD\rest S\deq\{X\cap S:\,X\in\DD\}.$$

\smallskip

We have fixed cardinals $\lambda,\mu$ and $\kappa$ at the beginning
of this section.
By [Sh 247, 6] or [Sh 212, 14], we know that 
$\DD_\lambda\rest S^\lambda_{\neq\kappa}$
is not $\lambda^+$-saturated.
We shall now introduce a sufficient condition,
under which we can prove that $\DD_\lambda\rest S^\lambda_\kappa$
is not $\lambda^+$-saturated. In Theorem 1.15, we extend the result to
a somewhat larger class of normal filters on $\lambda$.

\smallskip

{\bf Definition 1.1.} Suppose $\chi>\aleph_0$ is a regular cardinal and
$\eta>\chi$ is an ordinal. We say that $\eta$ {\it has the
strong non-reflection property for\/} $\chi$, if there is a
function $h:\eta\to\chi$ such that for every $\delta\in S^\eta_\chi$,
there is a club subset $C$ of $\delta$ with $h\rest C$ strictly increasing.
In such a case we say that $h$ witnesses the strong non-reflection
of $\eta$ for $\chi$.

If $\eta$ does not have the strong non-reflection property for $\chi$,
then we say that $\eta$ {\it has the weak reflection for\/} $\chi$,
or that $\eta$ {\it is weakly reflective for\/} $\chi$.

If $h:\eta\to\chi
$ is a function, we define
$$\ref(h)\deq\{\delta\in S^\eta_\chi: h\rest C\hbox{ is not strictly
increasing for any club } C \hbox{ of }\delta\}.$$

\smallskip

We shall first make some general remarks about Definition 1.1.

\smallskip

{\bf Observation 1.2}(1) If $\eta$ is weakly reflective for $\chi>
\aleph_0$,
and $\zeta>\eta$, {\it then\/} $\zeta$ is weakly reflective for $\chi$.

(2) If $\eta>\chi=\cf(\chi)>\aleph_0$, {\it then\/}

$\eta$ has the strong non-reflection
property for $\chi$ iff
there is an $h:\eta\to\chi$ such that for all $\delta\in
S^\eta_\chi,$
there is a club $ C$ of $\delta$ with
the property that $h\rest C$ is 1-1.

In fact, the two sides of the equivalence can be witnessed by the
same function $h$.

(3) Suppose that $\chi$ is a given uncountable
regular cardinal such that there is an
$\eta$ which weakly reflects at $\chi$.

{\it Then\/} the minimal such $\eta$, which we shall
denote by $\theta^\ast(\chi)$, is
a regular cardinal $>\chi$.

Consequently, $\zeta$ weakly reflects at $\chi$ iff
there is a regular
cardinal $\theta\leq\zeta$ such that $\theta$ weakly reflects at
$\chi$.

(4) Suppose that $\chi$ is a  regular cardinal, and 
$\eta=\theta^\ast(\chi)$ or
$\eta$ is an ordinal of cofinality $\ge\theta^\ast(\chi)$
such that $S^\eta_{\theta^\ast(\chi)}$
is stationary in $\eta$. (This makes sense, as by (3), if $\theta^\ast(\chi)$
is defined, then it is a regular cardinal.)

{\it Then\/}
not only that $\eta$ is weakly reflective for $\chi$,
but for every
$h:\eta\to\chi$,
the set $\ref(h)$ is stationary.

{\bf Proof.} (1) Follows from the definition.
 
(2) If $\eta$ has the strong non-reflection for $\chi$,
the other side of the above equivalence
is obviously true.

In the other direction,
suppose that $h$
satisfies the conditions on the right hand side,
and fix a $\delta\in S^\eta_\chi$.
Let $C$ be a club of $\delta$ on which $h$ is 1-1.
In particular note that ${\rm ran} (h)$ is cofinal in $\chi$ and that
$\otp(C)=\chi$.
By induction on $\gamma<\chi$,
define
$$\beta_\gamma=\Min\Bigl\{\alpha\in C:\,\bigl(\forall\beta\in(C\setminus\alpha)
\bigr)\,
\bigl(\forall\zeta<\gamma)\bigr)\,
\bigl(h(\beta)>h(\beta_\xi)\bigr)\Bigr\}.$$
Then $D=\{\beta_\gamma:\,\gamma<\chi\}$ is a club of $\delta$
and $h\rest D$ is strictly increasing.
 
(3) Let $\eta=\theta^\ast(\chi)$ be the minimal ordinal which weakly
reflects at $\chi$, so obviously $\eta>\chi$. 
If $\eta$ is not a regular cardinal, we can find
an increasing continuous sequence of ordinals
$\langle\eta_i:
\,i<\zeta\rangle$ which is cofinal in $\eta$, and such that 
$\zeta<\eta$. We can also assume that $\eta_0=0$ and $\eta_1>\chi$.
In addition,
for every $i<\zeta$, if $i$ is a successor, we can assume that $\eta_i$
is also a successor.

Then, for every $i<\zeta$, there is an $h_{i+1}
:\eta_{i+1}\to\chi$
which exemplifies
that $\eta_{i+1}$ has the strong non-reflection property for $\chi$.
There is also a $g:\zeta\to\chi$ which witnesses that $\zeta$
has the strong non-reflection property for $\chi$.
Define $h:\eta\to\chi$ by:
$$h(\alpha)=\cases{h_{i+1}(\alpha)   & if $\alpha\in(\eta_i,\eta_{i+1})$\cr
                   g(i)              & if $\alpha=\eta_i.$\cr}$$
Since $\eta$ is weakly reflective for $\chi$, there must
be a $\delta\in S^\eta_\chi$ such that for no club $C$ of
$\delta$, is $h\rest C$ strictly increasing. We can distinguish two cases:

\smallskip

$\underline{\hbox{Case 1.}}$ $\delta\in (\eta_i,\eta_{i+1}]$ for some
$i<\zeta$.

Then there is a club $C$ of $\delta$ on which $h_{i+1}$ is strictly increasing.
But then $h\rest (C\setminus
(\eta_i+1))=h_{i+1}\rest (C\setminus(\eta_i+1))$
is strictly increasing, and $C\setminus(\eta_i+1)$ is a club of
$\delta$. A contradiction.

$\underline{\hbox{Case 2.}}$ $\delta=\eta_i$ for some
limit $i<\zeta$.

Notice that 
$D=\{\eta_j:\,j<i\}$ is a club of $\eta_i$. We know that
there is a club $C$ of $i$ on which $g$ is increasing. Then
setting $E=C\cap D$, we conclude that
$h\rest E$ is strictly increasing.
Therefore, a contradiction.

(4) We first assume
that $\eta=\theta^\ast(\chi)$.
By (2),
$\theta^\ast(\chi)$ is necessarily a regular cardinal $>\chi$.
Let $h:\theta^\ast(\chi)\to\chi$ be given.
Suppose that $\ref(h)$ is non-stationary
in $\theta^\ast(\chi)$ and fix a club $E$ in $\theta^\ast(\chi)$ such that
$E\cap\ref(h)=\emptyset$. 
Without loss of generality, $\otp{E}=\theta^\ast(\chi)$.
Let us fix an increasing enumeration
$E=\{\alpha_i:\,i<\theta^\ast(\chi)\}$.
Without loss of generality, $\alpha_0=0$
and $\alpha_1>\chi$
and $\alpha_{i+1}$ is a successor for every $i$. So, for each $i$, there is a
function $h_i:
\alpha_{i+1}\to\chi$ which witnesses that 
$\alpha_{i+1}$ is strongly non-reflective at
$\chi$.

But now we can use $h$ and $h_i$ for $i<\theta^\ast(\chi)$
to define a function which will
contradict that $\theta^\ast(\chi)$ weakly reflects at $\chi$, similarly to the
proof of (3).

The other case is
that $\eta>\theta^\ast(\chi)$,
so $S^\eta_{\theta^\ast(\chi)}$ is stationary in $\eta$.
Let $h:\eta\to\chi$ be given.

If $\delta\in S^\eta_{\theta^{\ast}(\chi)}$, then let us fix a club $C_\delta$
of $\delta$ such that $\otp(C_\delta)=\theta^{\ast}(\chi)$.
Let $C_\delta=\{\alpha_i:\,i<\theta^{\ast}(\chi)\}$ be an increasing
enumeration. Then $h\rest C_\delta$ induces a function
$g_\delta:\theta^{\ast}(\chi)
\to\chi$, given by $g_\delta(i)=h(\alpha_i)$.

We wish to show that $\ref(h)$ is stationary in $\eta$. So, let $C$ be a club
of $\eta$. 
As we know that $S^\eta_{\theta^\ast(\chi)}$ is stationary in
$\eta$, we can find a $\delta\in S^\eta_{\theta^{\ast}(\chi)}$ which is an
accumulation
point of $C$. Then $C\cap C_\delta$ is a club of $\delta$, so
$E\deq\{i<\theta^{\ast}(\chi):\,
\alpha_i\in C\}$ is a club in $\theta^{\ast}(\chi)$.
Therefore, there is an $i\in E\cap\ref(g_\delta)$,
by the first part of this proof.
Then $\alpha_i\in\ref(h)\cap C$.\eop${}_{1.2.}$

\medskip

{\bf Remark 1.2.a.} One can ask the question of 1.2.3 in the opposite direction:
suppose
that $\sigma$ reflects at some $\eta$, what can we say about the first such
$\eta$?
This is an independence
question, for more on this see [CDSh 571].

\medskip

We 
find it convenient to introduce the following

\smallskip

{\bf Definition 1.3.} For ordinals $\eta>\chi>\aleph_0$,
where $\chi$ is a regular cardinal, we define
$$
\eqalign{\II[\eta,\chi)=\{A\subseteq\eta:\,&\hbox{there is a function }
h:\eta\to\chi\hbox{ 
such that for
every }\delta\in A\cap S^\eta_\chi,\cr
&\hbox{ there is a club } C\hbox{ of }
\delta
\hbox{ with } h\rest C\hbox{ strictly increasing.}\}\cr}$$
Saying that $\eta$ has the strong non-reflection property for $\chi$
is equivalent to claiming that $\eta\in\II[\eta,\chi)$,
or $S^\eta_\chi\in\II[\eta,\chi)$. In such a case, we say that
$\II[\eta,\chi)$ is {\it trivial\/}.

$I(\sigma,\chi)$ is the statement:

There is a $\eta\in(\chi,\sigma)$ with $\II[\eta,\chi)$ non-trivial.

\medskip

{\bf Theorem 1.4.} Suppose that $\eta>\chi>\aleph_0$ and $\chi$ is 
a regular cardinal.

{\it Then\/} $\II[\eta,\chi)$ is a $\chi$-complete
ideal on $\eta$.

{\bf Proof.} $\II[\eta,\chi)$ is obviously non-empty and downward
closed.

Suppose that $A_i\,(i<i^\ast<\chi)$ are sets from $\II[\eta,\chi)$,
and that $h_i:\eta\to\chi$ witnesses that $A_i\in\II[\eta,\chi)$
for $i<i^\ast$.  Let $A=\cup_{i< i^\ast}A_i$, and we shall see
that $A\in\II[\eta,\chi)$. This will be exemplified by the function
$h:\eta\to\chi$ defined by $h(\beta)\deq\sup\{h_i(\beta):\,i<i^\ast\}.$
Let $\delta\in A\cap S^\eta_\chi$ be given.

Then $\delta\in A_i$ for some $i<i^\ast$, and therefore there is a
club $C$ of $\delta$ such that $h_i\rest C$ is strictly increasing.
Then it must be that
$\otp(C)=\cf(\delta)=\chi$, and
we can enumerate $C$ increasingly
as $\{\beta_\epsilon:\,\epsilon<\chi\}.$
Therefore, the sequence $\langle h_i(\beta_\epsilon):\,
\epsilon <\chi\rangle$ is a strictly increasing sequence in $\chi$,
and $\chi=\sup_{\epsilon<\chi}\bigl(h_i(\beta_\epsilon)\bigr).$

On the other hand, for every $\epsilon$, we have that $h(\beta_\epsilon)
<\chi$, so there is a minimal $\xi(\epsilon)<\chi$ such that
$h(\beta_\epsilon)<h_i(\beta_{\xi(\epsilon)})$ for some $i<i^\ast$.
Let
$$E\deq\{\zeta<\rho:\,\bwedge_{\epsilon<\zeta}\xi(\epsilon)<\zeta\},$$
so $E$ is a club of $\chi$ and $C^\ast\deq\{\beta_\zeta:\,\zeta\in E\}$
is a club of $\delta$.

But then, for $\epsilon_1<\epsilon_2\in E$ we have
$\xi(\epsilon_1)<\epsilon_2$, so for some
$i<i^\ast$ we have $h(\beta_{\epsilon_1})<
h_i(\beta_{\xi(\epsilon_1)})<h_i(\beta_{\epsilon_2})
\leq h(\beta_{\epsilon_2})$, so 
$h\rest C^\ast$ is strictly
increasing.
\eop${}_{1.4.}$

\medskip

If there is a square on $\sigma^+$, then $I(\sigma^+,\chi)$ is false.
As there are various notations in use,
to make this statement precise, we state the definition of the $\square$
principle that we use. Note that what we refer to as $\square\langle
\sigma^+\rangle$,
some authors regard as $\square_\sigma$.

\smallskip

{\bf Definition 1.5.} Suppose that $\chi$
and $\sigma$ are cardinals. 
$\square_{(\chi,\sigma)}$ denotes the following
statement:

$\chi<\sigma$ and
there is a sequence $\langle C_\delta:\delta\in (SING
\cap LIM)
\,\&\,\chi<\delta<\sigma
\rangle$ such that

\item{(1)} $C_\delta$ is a club in $\delta$,
\item{(2)} $\otp(C_\delta)<\delta$,
\item{(3)} If $\delta$ is an
accumulation point of $C_\alpha$, then $C_\delta$ is defined
and $C_\delta=\delta\cap C_\alpha$. In particular, $\delta>\chi.$

We use
$\square
\langle\chi^+\rangle$ as a shorthand for $\square_{(\chi,\chi^+)}.$
The sequence as above is called a $\square_{(\chi,\sigma)}$ sequence,
and its subsequences are called partial
$\square_{(\chi,\sigma)}$-sequences. 

\medskip

{\bf Observation 1.6}(0) Note that by a closed unbounded set of $\omega$,
we simply mean an unbounded subset of $\omega$. Similarly,
any $\delta$ with $\cf(\delta)=\aleph_0$ will have a club subset consisting
of an unbounded $\omega$-sequence in $\delta$.
So,
$\square\langle\omega_1\rangle$ trivially holds.

(1) If $\chi_1\leq\chi_2<\sigma$,
{\it then\/}  
$\square_{(\chi_1,\sigma)}\implies \square_{(\chi_2,\sigma)}.$
If $\chi<\sigma_1\le\sigma_2$,
{\it then\/}
$\square_{(\chi,\sigma_2)}\implies
\square_{(\chi,\sigma_1)}$.


\medskip

{\bf Theorem 1.7.} If
$\theta>\chi$
is regular 
, {\it then\/} $\square_{(\chi,\theta)}$ implies that
$\theta$ has the strong
non-reflection property for $\chi$.

{\bf Proof.}
\relax From [Sh -g VII 1.7], we recall the following Fact 1.7.a. For the reader's
convenience, we also include the proof.

\smallskip

{\bf Fact 1.7.a.} Assume that $\square_{(\chi,\theta)}$ holds.

{\it Then,\/} there is a partial $\square_{(\chi,\theta)}$-sequence 
$$\langle C_\delta:\,\delta\in (SING\cap
LIM)\,\&\,\chi<\delta<\theta\,\&\,\cf(\delta)<\chi\rangle,$$
such that for each $\delta$ for which $C_\delta$ is defined,

\item{$(i)$} $C_\delta$ is a club subset of $\delta$.
\item{$(ii)$} $\otp(C_\delta)<\delta\quad\&\quad \cf (\delta)<\chi
\implies\otp (C_\delta)<\chi$.
\item{$(iii)$} $\gamma\in\acc(C_\delta)
\quad\&\quad C_\gamma\hbox{ defined }\implies C_\gamma=C_\delta\cap
\gamma.$

{\bf Proof of 1.7.a.} We 
start with a sequence $\langle D_\delta:\,\delta\in (SING\cap LIM)
\,\&\,\chi<\delta<\theta\rangle$ which exemplifies $\square_{(\chi,\theta)}
$.

We can without loss of generality assume that each $D_\delta$ satisfies
$D_\delta\cap\chi=\emptyset.$
For each $\delta$ for which $D_\delta$ is defined, we can define
a 1-1 onto function $f_\delta:\,D_\delta\to\otp(D_\delta)$
by
$$f_\delta(\alpha)\deq\otp(D_\delta\cap\alpha).$$
Note that $\gamma\in\acc(D_\delta)\implies f_\delta\rest D_\gamma=f_\gamma.$
Now we define $C_\delta$ for $\delta<\theta$ with $\cf(\delta)<\chi$
by induction on $\delta<\theta$.

If $D_\delta$ is not defined, then $C_\delta$ is not either.

If $D_\delta$ is defined, and $D_{\otp(D_\delta)}$ is not defined, 
we set $C_\delta=D_\delta$.
Note that $\otp(D_\delta)<\chi$ in this case.

Finally, suppose that
both $D_\delta$ and $D_{\otp(D_\delta)}$ are defined. Then
$\delta>\otp(D_\delta)>\chi$ and
$\cf\bigl ((\otp(D_\delta)\bigr)=\cf(\delta)<\chi$. So
$C_{\otp(D_\delta)}$ is already defined and
we can define
$$C_\delta\deq\{\alpha\in D_\delta:\,f_\delta(\alpha)\in C_{\otp(D_\delta)}\}.$$
We can check that
$$\bar{C}\deq\langle C_\delta:\,\delta\in (SING\cap LIM)\,\&\,\chi
<\delta<\theta\,\&\,cf(\delta)<\chi\rangle$$
is as required.
One thing to note is that if
$\gamma\in\acc(C_\delta)$ and $D_{\otp(D_\delta)}$
is defined, then $D_{\otp(D_\gamma)}$ must be
defined too.\eop${}_{1.7.a.}$

\smallskip

Fixing a sequence $\bar{C}$ like in Fact 1.7.a,
the following defines a function $h:\theta\to\chi$:
$$h(\delta)=\cases{\otp(C_\delta) & if $\delta\in S^\theta_{<\chi}
\setminus\chi$\cr
                                  0 & otherwise.\cr}$$

Now, if $\delta\in S^\theta_{\chi}$,
we can choose a club $E_\delta$ of $\delta$ which consists only
of elements of $S^\theta_{<\chi}\setminus\chi.$ We let
$$D_\delta=\acc(C_\delta)\cap E_\delta.$$
Then $h\rest D_\delta$ is increasing and $D_\delta$ is a club in
$\delta$.\eop${}_{1.7.}$

\medskip

We can also show,
for example, that $I(\sigma,\chi^+)$ is consistently true, if
$\sigma>\chi^+$ and $\chi$ is regular. This follows from
Fact
1.10 below. This Fact also explains why we choose
to call the properties under consideration
``weak reflection'' and ``strong non-reflection''. 
There are other reflection properties that imply the
weak reflection, like the reflections considered in [Ba] and elsewhere,
so Fact
1.10 can be used as an example
of a proof that the weak reflection is weaker than
other reflection principles. Let us first recall
the notion of stationary reflection.

\smallskip

{\bf Definition 1.8}(0) Suppose that $S$ is a stationary subset of 
a cardinal $\theta$ with $\cf(\theta)>\aleph_0$.
We say that $S$ {\it reflects\/} at $\delta\in\theta$, if
$\cf(\delta)>\aleph_0$ and $S\cap\delta$
is stationary in $\delta$. 
We say that $S$ {\it is reflecting\/} if there is $\delta\in \theta$
such that $S$
reflects at $\delta$. Otherwise, $S$ is said to be
{\it non-reflecting\/}.

(1) A regular cardinal $\theta$ is {\it reflecting\/} iff for every regular
$\chi$ such that $\chi<\chi^+<\theta$, every stationary $S\subseteq
S^\theta_\chi$ is reflecting. 

(2) For a regular cardinal $\theta$, 
notation $REF(\theta)$ means that $\theta$ is
reflecting. $REF$ is the statement denoting that for every 
$\theta>\aleph_1$ which is a regular cardinal, $REF(\theta)$ holds.

(3) Suppose that $\lambda>\theta,\kappa$ are regular cardinals
and $\kappa>\aleph_0$. We define
the statement $\hbox{Ref}(\lambda,\kappa,\theta)$ to mean:

For every $S\subseteq S^\lambda_\theta$ which is stationary,
there is a $\delta$ of cofinality $\kappa$, such that $S
$ reflects at $\delta$.

\smallskip

{\bf Remark 1.9} (0) If
$\delta$ is an ordinal of uncountable cofinality,
then $\delta$ has a club subset consisting only of elements of
cofinality $<\cf(\delta)$. Therefore, if $S$ reflects at $\delta$,
then $S$ has to have elements of cofinality $<\cf(\delta)$.
This explains the gap of one between $\chi$ and $\theta$ in the definition
of $REF(\theta)$.

\smallskip

$REF$ is consistent modulo the existence of infinitely
many supercompact cardinals [Sh 351]. 
The consistency of $\hbox{Ref}(\lambda,\kappa,\theta)$
has also been extensively studied,
starting with a
result of
J. Baumgartner in [Ba] 
that $CON(\hbox{Ref}(\aleph_2,\aleph_1,\aleph_0))$
follows from the existence of a weakly compact cardinal.

(1) M. Magidor points out the following equivalent definition of the
weak reflection, from which it is easy to see that it is weaker than
what is usually meant by reflection:

We can say that $\theta>\kappa$ weakly reflects at $\kappa$
iff for any partition
$$S^\theta_{<\kappa}=\cup_{i<\kappa}S_i,$$
there is an $i<\kappa$ and an
$ \alpha\in S^\theta_\kappa$
such that $S_i\cap\alpha$ is stationary in $\alpha$.


\medskip

{\bf Fact 1.10}(0) Suppose that
$\chi$
is a successor cardinal.

{\it Then\/}
$$REF(\chi^+)\implies\chi^+\hbox{ weakly reflects at }\chi.$$

(1) $\hbox{Ref}(\lambda,\kappa,\theta)\implies\lambda\notin\II[\lambda,\kappa)$.

{\bf Proof.} (0)
Let $\chi=\sigma^+$.
So,
suppose
for contradiction that $h:\chi^+\to\chi$ is a function such that for
every $\delta\in S^{\chi^+}_\chi$, there is a club of $C_\delta$
with $h\rest C_\delta$ is 1-1.

Now, $h$ is regressive on $(\chi,\chi^+)$.
Since $(\chi,\chi^+)\cap S^{\chi^+}_\sigma$
is stationary in $\chi^+$, there is a stationary $S\subseteq (\chi,\chi^+)
\cap S^{\chi^+}_\sigma$ such that $h\rest S$ is a constant.

By $REF(\chi^+)$, there is a $\delta\in\chi^+$ such that $S\cap\delta$
is stationary in $\delta$. As in 1.9, we conclude that
$\cf(\delta)>\sigma.$ Therefore $\cf(\delta)=\chi$ and $C_\delta$
is defined. 
But then $S\cap C_\delta$
is stationary in $\delta$, and $h\rest (S\cap C_\delta)$ is constant.

(1) Similar.
\eop${}_{1.10.}$  

\medskip

We now go on to present the last two facts before
we proceed to the Main Theorem.

\medskip

{\bf Fact 1.11.} Suppose that $\cf(\delta)>\aleph_0$ and $f
:\delta\to\delta$ is a function
which is
not increasing on any club of $\delta$. {\it Then\/}
there is a stationary set $S$ in
$\delta$ such that $\alpha\in S\implies f(\alpha)<
\Min(S)$. (Hence, there
is also a stationary subset of $\delta$ on which $f$ is a constant.)

{\bf Proof.} We fix an increasing continuous sequence of ordinals
$\langle \alpha_\epsilon:\,\epsilon<\cf(\delta)\rangle$ which is
cofinal in $\delta$. Let $T\deq\{\epsilon:\,f(\alpha_\epsilon)<\alpha_\epsilon
\}$. So $T\subseteq\cf(\delta).$

We shall see that $T$
is stationary in $\cf(\delta)$.
Let us first assume that it is true, and
define for $\epsilon\in T\setminus\{0\}$,
$$g(\epsilon)\deq\Min\{\zeta:\,f(\alpha_\epsilon)<\alpha_{\zeta+1}\}.$$
Therefore, $g$ is regressive, and we can find a stationary $T_1\subseteq T$
such that $g\rest T_1$ is constantly equal to some $\zeta^\ast$.

Let $S\deq\{\alpha_\epsilon:\,\epsilon\in T_1\setminus(\zeta^\ast +2)\}.$
We can check that this $S$ is as required.

It remains to be seen that
$T$ is stationary in $\cf(\delta)$. Suppose not.
Then we can find a club $C$ in $\cf(\delta)$ such that $C\cap T=\emptyset$.
Let
$$E\deq\{\epsilon:\,\epsilon\in LIM\cap C\,\&\,(\forall\zeta<\epsilon)(
f(\alpha_\zeta)<\alpha_\epsilon)\}.$$
Then $E$ is also a club of $\cf(\delta)$. But then
$D\deq\{\alpha_\epsilon:\,\epsilon\in E\}$ is a club of $\delta$ and
$f\rest D$ is strictly increasing, contradicting our assumption.\eop${}_{1.11.}$

\medskip

{\bf Remark 1.12.} 
As a remark on the side:
we cannot improve the previous result to conclude, from
the assumptions given above, that there is a club $C$
of $\delta$ such that $f\rest C$ is
constant. Namely, if we take a stationary costationary set $S$ in $\delta$,
and define $f$ on $\delta$ by
$$f(\alpha)=\cases{1 & if $\alpha\in S$\cr
                   0 & otherwise,\cr}$$
then this $f$ is neither increasing nor constant on any club of $\delta$.
In the following Fact 1.12.a we recall
a more general result along the lines of Fact 1.11.
This Fact is not used in the proof of the Main Theorem 1.13,
and a reader who is in a hurry may
without loss of continuity proceed directly to 1.13.

\medskip

{\bf Fact 1.12.a.} Suppose that $\delta$ is an ordinal with
$\cf(\delta)>\aleph_0$ and $f$ is a function from $\delta$ to the
ordinals. {\it Then,\/}
there is a stationary $S\subseteq \delta$ such that

{\it either\/} $f\rest S$ is constant,

{\it or\/} $f\rest S$ is strictly increasing.

{\bf Proof.} Let $\theta=\cf(\delta)>\aleph_0$, and
$\langle \alpha_\epsilon:\,\epsilon<\theta\rangle$ a strictly
increasing 
enumeration of a club of $\delta$.
Let $E\deq\theta\cap LIM$.
We define a partial function $g:E\to\theta$ as follows:
$$g(\epsilon)=\zeta \hbox{ if }\cases{ &(a) $f(\alpha_\epsilon)\leq
                                        f(\alpha_\zeta)$\cr
                              &    (b) $\bigl(\forall\xi<\epsilon\bigr)
                            \bigl(f(\alpha_\epsilon)\le f(\alpha_\xi)
                            \implies f(\alpha_\zeta)\le f(\alpha_\xi)\bigr)$\cr
                              &    (c) $\zeta<\epsilon$
is minimal under (a) and (b).\cr}$$
Note that $g(\epsilon)<\epsilon$ for all $\epsilon\in{\Dom}(g)$.
Now we consider three cases:

\smallskip

$\underline{\hbox{Case 1.}}$ $S_0\deq E\setminus{\Dom}(g)$ is stationary
in $\theta$.

Then $S\deq\{\alpha_\epsilon:\,\epsilon\in S_0\}$
is stationary in $\delta$. We claim that
$f\rest S$ is strictly
increasing on $S$. Otherwise, there would be an $\epsilon\in S$
such that there is a $\zeta_0<\epsilon$ with $f(\alpha_\epsilon)\le
f(\alpha_{\zeta_0}),$
which contradicts the fact that $\epsilon\notin
{\Dom}(g).$

$\underline{\hbox{Case 2.}}$ ${\Dom}(g)$ is stationary in $\theta$
and
$S_1\deq\{\epsilon\in{\Dom}(g):\,f(\alpha_\epsilon)=f(\alpha_{g(\epsilon)})
\}$ is stationary in $\theta$.

Since $g\rest S_1$ is regressive, there is a stationary $S_2
\subseteq S_1$ such that $g\rest S_2$ is constant. Then
$S\deq\{\alpha_\epsilon:\,\epsilon\in S_2\}$ is stationary in $\delta$,
and $f\rest S$ is constant.

$\underline{\hbox{Case 3.}}$ ${\Dom}(g)$ is stationary in $\theta$,
but $S_1$ is not stationary.

Then $S_3\deq{\Dom}(g)\setminus S_1=\{\epsilon\in {\Dom}(g):\,
f(\alpha_\epsilon)<f(\alpha_{g(\epsilon)})\}$ is stationary,
and there is a stationary $S_4\subseteq S_3$ such that $g\rest S_4$
is a constant. Let $S=\{\alpha_\epsilon:\,\epsilon\in S_4\}$,
so $S$ is a stationary subset of $\delta$. We claim that
$f\rest S$ is strictly increasing.

Otherwise, there are $\epsilon_1<\epsilon_2\in S_4$ such that
$f(\alpha_{\epsilon_2})\leq f(\alpha_{\epsilon_1}).$
On the other hand, $f(\alpha_{\epsilon_1})<
f(\alpha_{g(\epsilon_1)})=f(\alpha_{g(\epsilon_2)})$, since
both $\epsilon_1$ and $\epsilon_2$ are members of $S_4$.
This contradicts the definition of $g(\epsilon_2)$, since $\epsilon_1<
\epsilon_2$.\eop${}_{1.12.a.}$

\medskip

We now present our main result.

\medskip

{\bf Main Theorem 1.13.} Assume that $\lambda=\mu^+$ and $\mu>\cf(\mu)
=\kappa >\aleph_0$.
In addition, for some $\theta\in(\kappa,\lambda)$,
we know that $\theta$ has the
weak reflection property for $\kappa$.

{\it Then\/} $\DD_\lambda\rest S^\lambda_\kappa$ is not $\lambda^+$-saturated.

{\bf Proof.} By 1.2.3, without loss of generality, $\theta$ is a regular 
cardinal, so $\theta<\mu$.

\smallskip

{\bf Fact 1.13.a.} There is a stationary subset $S$ of $S^\lambda_\theta$
such that:

For some $S^+\supseteq S$ and $S^+\subseteq\lambda\setminus\kappa$, 
there is a sequence
$$\bar{C}=\langle C_\alpha:\alpha\in S^+\rangle,$$ such that, for every
$\alpha\in
S^+$,

\item{1.} $C_\alpha$ is a subset of $\alpha
\setminus\kappa$ and $\otp(C_\alpha)\leq\theta$.
\item{2.} $\alpha$ is a limit ordinal $\implies\sup(C_\alpha)=\alpha$.
\item{3.} $\beta\in C_\alpha\implies(\beta\in S^+$ and $C_\beta=C_\alpha\cap
\beta)$.
\item{4.} $\cf(\alpha)=\theta$ iff $\alpha\in S$.
\item{5.} For every club $E$ in $\lambda$, there are stationarily many
$\delta\in S,$ such that for all $\alpha,\beta$:
$$(\alpha<\beta\,\,\&\,\,
\alpha,\beta\in C_\delta)\implies(\alpha,\beta]\cap E\neq\emptyset.$$

{\bf Proof of 1.13.a.} 
Since $\theta^+<\lambda$ we can apply
[Sh 420, 1.5], so
there is a stationary subset
$S_1$
of $S^\lambda_\theta\setminus\{\theta\}$ with $S_1\in I[\lambda]$.
By [Sh 420, 1.2] this means that
there is a sequence $\langle D_\alpha:\,\alpha< \lambda
\rangle$
such that

\item{(a)} $D_\alpha$ is a closed subset of $\alpha$.
\item{(b)} $\alpha^\ast\in\nacc(D_\alpha)\implies 
D_{\alpha^\ast}
=D_\alpha\cap\alpha^\ast$.
\item{(c)} For some club $E$ of $\lambda$, for every $\delta\in S_1\cap E$,
$$\delta=\sup(D_\delta)\quad\&\quad
\otp(D_\delta)=\theta.$$
\item{(d)} $\hbox{nacc}(D_\alpha)$ is a set of successor ordinals.

\smallskip

{\it Observation.\/}
We do not lose generality if we in addition require that for each $\alpha$,
$\otp(D_\alpha)\leq\theta$.

[Why? Let $E$ be the club guaranteed by (c). 
Since $S_1$ is stationary, so is $S_1\cap E$, so we can define
for $\alpha\in\lambda$
$$D_\alpha^\dagger=\cases{D_\alpha  & if $\alpha\in S_1\cap E$\cr
                                    & or $\exists\beta>\alpha\,\bigl
                                        (\beta\in S_1
                                            \cap E\,\&\, \alpha\in 
                                      \nacc(D_\beta)\bigr).$\cr
                         \emptyset  & otherwise.\cr}$$
Then we can set $S_2\deq S_1\cap E$, 
so the sequence
$\langle D_\alpha^\dagger:\,\alpha<\lambda\rangle$ will satisfy (a)--(d),
with $S_1$ replaced by $S_2$, and $\otp(D_\alpha^\dagger
)\le\theta$ will hold
for each $\alpha$.]

{\it Continuation of the Proof of 1.13.a.\/}
Now, for any club $F$ of $\lambda$, we let
$$\bar{C}[F]\deq\langle C_\delta[F]:\,\delta\in
S^+[F]\deq S_1\cap E\cap\acc(F)\cup(\hbox{successors})\setminus
\kappa\rangle,$$
where
$$C_\delta[F]\deq\Bigl\{\gamma\in\nacc (D_\delta)
\setminus\kappa:\,F\cap\bigl[\sup(\gamma
\cap D_\delta),\gamma\bigr)\neq\emptyset\Bigr\}.$$
We claim that $\bar{C}$ can be set to be equal to $\bar{C}(F)$ for some 
club $F$, with $S=S[F]\deq E\cap S_1\cap\acc(F)
\setminus\kappa$ and $S^+=S^+[F]$.  
We are using the ideas of [Sh 365,\S2].

It is easily checked that (1)--(4) are satisfied for any club $F$.
So, let us suppose that (5) is not satisfied for any choice of $F$,
and we shall obtain a contradiction.

By induction on $\zeta<\theta^+$, we define a club $F_\zeta$
of $\lambda$.

If $\zeta=0$, we let $F_\zeta=LIM$.

If $\zeta$ is a limit ordinal, we let $F_\zeta=\cap_{\xi<\zeta}F_\xi.$
This is still a club of $\lambda$, as $\zeta<\theta^+<\lambda$.

If $\zeta=\xi+1$, we have assumed that there is a club $F_\zeta$
such that the set
$$G_\xi[F_\zeta]\deq\{\delta\in S[F_\zeta]
:\,\forall \alpha ,\beta\in C_\delta[F_\xi]\,
\bigl(\alpha<\beta\implies (\alpha,\beta]\cap F_\zeta\neq\emptyset
\bigr)\}$$
is non-stationary. Without loss of generality, we assume that
$F_\zeta\subseteq F_\xi.$

At the end, let us let $C=\cap_{\zeta<\theta^+} F_\zeta.$ This is still a
club of $\lambda$, as $\theta^+<\lambda$. Since $S\deq S[C]$ is
stationary in $\lambda$, and for every $\zeta<\theta^+$, we have that
$S\subseteq S[F_\zeta]$, we conclude that the set
$$T\deq C\cap S\setminus (\cup_{\zeta<\theta^+} G_\zeta[F_{\zeta+1}])$$
is stationary. So, let us take a $\delta\in T\cap LIM$. Then $D_\delta$
is a club of $\delta$, since $T\cap LIM\subseteq S_1\cap E$.

For $\beta\in D_\delta$, we consider the sequence
$$\langle\sup(\beta\cap F_\zeta):\,\zeta<\theta^+\rangle.$$
This is a non-increasing sequence of ordinals $\le\beta$, so there must be a
$\zeta_\beta<\theta^+$ and $\gamma_\beta\le\beta$ such that
$$\zeta_\beta\leq\zeta<\theta^+\implies\sup(\beta\cap F_\zeta)
=\gamma_\beta.$$

Notice that $\zeta^\ast\deq\sup\{\zeta_\beta:\,\beta\in D_\delta\}
<\theta^+$, since $\otp(D_\delta)=\theta$.

Since $\delta\notin G_{\zeta^\ast}[F_{\zeta^\ast+1}]$,
by the definition of $G_{\zeta^\ast}[F_{\zeta^\ast+1}]$, there are
$\alpha<\beta\in C_\delta[F_{\zeta^\ast}]$ such that
$$(\alpha,\beta]\cap F_{\zeta^\ast+1}=\emptyset.$$
On the other hand, $\beta\in C_\delta[F_{\zeta^\ast}]\implies
\bigl[\sup(\beta\cap D_\delta),\beta\bigr)\cap F_{\zeta^\ast}\neq\emptyset
\implies\sup(\beta\cap D_\delta)\leq \sup(\beta\cap F_{\zeta^\ast})
=\gamma_\beta=\sup(\beta\cap F_{\zeta^\ast +1}).$ But
$\alpha<\beta$ and $\alpha\in\nacc(D_\delta)$, so
$\sup(\beta\cap D_\delta)\ge\alpha$. Therefore $\alpha\le\gamma_\beta
\le\beta$.
Note now that $\gamma_\beta$ must be a limit ordinal, since
$F_{\zeta^\ast+1}\subseteq LIM$. 
But $\alpha$ is a successor,
by (d) in the definition
of $D_\delta$. So, $\alpha<\gamma_\beta\le\beta\in F_{\zeta^\ast+1}$ and
$(\alpha,\beta]\cap F_{\zeta^\ast +1}=\emptyset,$ by
the choice of $F_{\zeta^\ast}$. A contradiction.\eop${}_{1.13.a.}$

\smallskip

{\it Continuation of the proof of 1.13.\/}
Let us fix $S, S^+$ and $\bar{C}\deq\langle C_\alpha:\alpha\in S^+
\rangle$ as in Fact 1.13.a. 
Denote by
$\cl (\bar{C})$ the following
$$\cl(\bar{C})=\{ C:\, C\subseteq S^+\,\wedge\,\forall\,\beta\in
 C\,(C_\beta=C\cap\beta)\}.$$
Now fix the following enumerations:

For $\alpha\in S^+$, let 
$$C_\alpha=\{\gamma(\alpha,\epsilon):\epsilon<\otp(C_\alpha)\},$$
such that $\gamma(\alpha,\epsilon)$ is increasing in $\epsilon$. For each
$\epsilon$, let 
$$\gamma^\ast(\alpha,\epsilon)\deq\cup_{\xi\leq
\epsilon}\gamma(\alpha,\xi).$$

Let $\mu=\Sigma_{i<\kappa}\mu_i$ 
be such that $\langle\mu_i:\, i\in\kappa\rangle$
is a continuous increasing sequence of cardinals, and,
for simplicity, $\mu_0>\theta$.

\smallskip

{\bf Claim 1.13.b.} Given enumerations as above,
we can for each $\alpha<\lambda$, find sets $a^\alpha_i (i<\kappa)$
such that 

\item{(I)} $\alpha=\cup_{i<\kappa}a^\alpha_i$.
\item{(II)} $\card{a^\alpha_i} 
\leq \mu_i$.
\item{(III)} $i<j\implies a^\alpha_i\subseteq a^\alpha_j $
and for $i$ a limit ordinal $<\kappa$,
$$a_i^\alpha=\cup_{j<i} a_j^\alpha.$$
\item{(IV)} $\beta\in a^\alpha_i\implies a^\beta_i\subseteq a^\alpha_i.$
\item{(V)} $\beta\in a^\alpha_i\cap S^+\implies C_\beta\subseteq a^\alpha_i$.

{\bf Proof of Claim 1.13.b.} It is easy to see that we can choose $a^\alpha_i$
for $\alpha<\lambda$ and $i<\kappa$ such that (I)--(III) are satisfied.
Suppose we have done so.
Then define
by induction on $\alpha<\lambda$, and then by induction on $j<\kappa$,
sets
$${a^\alpha_i}^\dagger \deq a_i^\alpha\bigcup
\bigl\{{a^\beta_i}^\dagger:\,\beta\in
                                         a^\alpha_i\cup\cup\{
                                         C_\beta:\,\beta\in a^\alpha_i
                                         \cap S^+\}\bigr\}
\bigcup\cup\bigl\{ C_\beta:\,\beta\in a^\alpha_i\cap S^+\bigr\}.$$
Now we can check that, 
by renaming $a^\alpha_i={a^\alpha_i}^\dagger$, we
satisfy the Claim.\eop${}_{1.13.b.}$

\smallskip

{\it Continuation of the proof of 1.13.\/}
So we now fix $a^\alpha_i$ as in Claim 1.13.b.

Similarly,
for each $\xi<\lambda^+$, we can let $\xi=\cup_{\alpha<\lambda}b^\xi_\alpha$,
where $b^\xi_\alpha$ are
$\subseteq$--increasing in $\alpha$,
while $\card{b^\xi_\alpha}
<\lambda$ and, if $\zeta\in b^\xi_\alpha$, then $b^\zeta_\alpha\subseteq
b^\xi_\alpha$. We
require in addition that, if $\xi=\zeta+1$, then $\zeta\in b^\xi_0$, and
if $\cf(\xi)<\lambda$, then $\xi=\sup(b^\xi_0)$.

We now define, for $\xi<\lambda^+$, a function $h_\xi:\lambda\to\lambda$. 
This function is given by
$h_\xi(\alpha)\deq\otp (b^\xi_\alpha)$.

Clearly, each $h_\xi$ is non-decreasing, and
$$\zeta\in b^\xi_\alpha\implies h_\zeta(\alpha)
< h_\xi(\alpha). \eqno{(a)}$$  

Now, fix a
sequence $\bar{E}=\{E_\eta:\,\eta\in S^\lambda_\kappa\}$,
where each $E_\eta
$ is a club subset of $\eta$ with 
$\otp(E_\eta)=\kappa$ and consisting only of elements of cofinality $<\kappa$.
We prove the following claim, with the intention of applying it later
on the functions $h_\xi\,(\xi<\lambda^+).$

\smallskip
 
{\bf Claim 1.13.c.} Suppose that $h:\lambda\to\lambda$ is non-decreasing,
and $\bar{E}$ is given as above.

{\it Then \/},
there is an $i=i(h)<\kappa$ such that for stationarily many
$\eta\in S_\kappa^\lambda$,
there is a $C=C(\eta)\in\cl({\bar{C}})$ with:

\item{(i)} $C\subseteq\eta=\sup(C).$
\item{(ii)} $C\subseteq\cup_{\beta\in E_\eta} a^\beta_i$ (hence
$\cup_{\gamma\in C}a^\gamma_i\subseteq \cup_{\beta\in E_\eta} a^\beta_i$,
by the choice of $a$'s).

\item{(iii)} $\eta=\sup \{\alpha\in C: h(\alpha)\in \cup_{\beta\in E_\eta}
 a^\beta_i$ and $h(\sup (C\cap\alpha))\in \cup_{\beta\in E_\eta} a^\beta_i\}$.

{\bf Proof of Claim 1.13.c.} Recall
the definition of $S$ and $C_\delta$'s from
Claim 1.13.a. In particular, for all $\delta\in S\subseteq 
S^\lambda_\theta$,
$$C_\delta=\{\gamma(\delta,\epsilon):\,\epsilon <\otp(C_\delta)=\theta\}.$$
Let us first fix a $\delta\in S$ and suppose that
for all
$\alpha\in C_\delta$, 
$$h(\alpha)<\Min (C_\delta\setminus(\alpha+1)).\eqno{(b)}$$
We define a function
$f=f_\delta: \theta\to\kappa$ by:

If $\epsilon\in(\theta\setminus S^\theta_{<\kappa})$, $f(\epsilon)=0$.

If $\epsilon\in S^\theta_{<\kappa}$,
$$f(\epsilon)=\Min\{ i<\kappa : h(\gamma^\ast(\delta,\epsilon)), h(\gamma(\delta,
\epsilon))\in a_i^{\gamma(\delta,\epsilon +1)}$$
$$\hbox{      and }
a_i^{\gamma^\ast(\delta,\epsilon)}\cap C_\delta\hbox{ is unbounded in }
\gamma^\ast(\delta,\epsilon)\}.$$

Remember that all $\epsilon\in S^\theta_{<\kappa}$
satisfy
$\aleph_0\leq\cf(\gamma^\ast(\delta,\epsilon))<
\kappa$.
By the
definition of $\gamma^\ast(\delta,\epsilon),$
we have that $C_\delta\cap\gamma^\ast(\delta,\epsilon)$ is
unbounded in $\gamma^\ast(\delta,\epsilon)$. Then, as $\langle 
a_i^{\gamma^\ast(\delta,\epsilon)}:\, i<\kappa\rangle$ is a 
$\subseteq$--increasing
sequence of sets with union
$\gamma^\ast(\delta,\epsilon)$,
it must be that $C_\delta\cap a_i^{\gamma^{\ast}
(\delta,\epsilon)}$ is eventually unbounded in $\gamma^{\ast}
(\delta,\epsilon)$.

On the other hand, since $h$ is non-decreasing and
$\gamma^{\ast}(\delta,\epsilon)\leq\gamma(\delta,\epsilon)
<\gamma(\delta, \epsilon +1)$,
we have  by $(b)$ (at the beginning of this proof) that
$h(\gamma^\ast(\delta,\epsilon)), h(\gamma(\delta,
\epsilon))<\gamma(\delta,\epsilon+1),$ so $h(\gamma^\ast(\delta,\epsilon)),
h(\gamma(\delta,\epsilon))
\in a_i^{\gamma(\delta,\epsilon +1)}$ for every large enough $i<\kappa$.

So, $f(\epsilon)$ is well defined.

\smallskip

We have assumed that $\theta$ is weakly reflective for $\kappa$. So,
by Observation 1.2.4, for  stationarily many
$\epsilon\in S^\theta_\kappa$,
$f$ is not strictly increasing on any club of $\gamma^\ast(\delta,\epsilon)$.
If we take any such $\epsilon$,
then $\cf\bigl(\gamma^\ast(\delta,
\epsilon)\bigr)>\aleph_0$
and $f\rest\gamma^\ast(\delta,\epsilon):
\,\gamma^\ast(\delta,\epsilon)\to\gamma^\ast(\delta,\epsilon)$,
as $\gamma^\ast(\delta,\epsilon)>\kappa$. So Fact 1.11 applies and
$f$ is constant on some
stationary subset 
of $\gamma^\ast(\delta,\epsilon)$. Let us denote that constant
value by $i_\epsilon=i_\delta(\epsilon)$.
Let $\eta\deq\gamma^\ast(\delta,\epsilon)$.

{\it Observation.\/}
The set
$\cup_{\alpha\in E_\eta} a^\alpha_{i_\epsilon}$ is an unbounded
(even stationary)
subset of $C_\delta\cap\gamma^\ast(\delta,\epsilon)$. 

[Why? We can check that
$e_\epsilon\deq\{\zeta<\epsilon:\,\gamma^\ast(\delta,\zeta)\in E_\eta\}$
is a club of $\epsilon$, so $s_\epsilon\deq\bigl\{\zeta\in e_\epsilon:\,f
(\zeta)=i_\epsilon\bigr\}\subseteq e_\epsilon$ is stationary,
and we conclude
that $\{\gamma(\delta,\zeta):\,\zeta\in s_\epsilon\}\subseteq\gamma^{\ast}
(\delta,\epsilon)$ is stationary in $\gamma^{\ast}(\delta,\epsilon).$
Note that for each $\zeta\in e_\epsilon$, by the choice of $\bar{E}$,
we have that $\cf(\gamma^\ast(\delta,\zeta))<\kappa$.

Now we take any $\zeta\in s_\epsilon$.
Since $\cf(\gamma^\ast(\delta,\zeta))<\kappa$, by the definition of $f$,
we have that $a_{i_\epsilon}^{\gamma^\ast(\delta,\zeta)}\cap C_\delta$
is unbounded in $\gamma^\ast(\delta,\zeta)$. Then
$$\cup_{\alpha\in E_\eta} a_{i_\epsilon}^\alpha\cap C_\delta\supseteq
\cup_{\gamma^\ast(\delta,\zeta)\in s_\epsilon}
a_{i_\epsilon}^{\gamma^\ast(\delta,\zeta)}
\cap C_\delta,$$
and since $s_\epsilon$
is stationary in $\gamma^\ast(\delta,\epsilon)$, we derive
the desired conclusion.]

{\it Continuation of the Proof of 1.13.c.\/}
So, by (V) in the choice of $a$'s and (3) of Fact 1.13.a,   
we conclude that
$\cup_{\alpha\in E_\eta} a^\alpha_{i_\epsilon}$ contains the entire
$C_\delta\cap\gamma^\ast(\delta,\epsilon)$. Hence, 
by (IV) of 1.13.b,
$$\cup_{\alpha\in E_\eta} a^\alpha_{i_\epsilon}\supseteq
\cup\{ a^\beta_{i_\epsilon}:\beta\in C_\delta\cap\gamma^\ast(\delta,\epsilon)
\}.$$

Now, there must be some $j=j_\delta<\kappa$,
such that $\{\epsilon\in S^\theta_\kappa
:\, i_\epsilon\hbox{ is well
defined and }=j\}$ is a stationary subset of $\theta$.

If we let $\delta$ vary,
then
$j_\delta$ is defined for every $\delta\in S^\lambda_\theta$ which satisfies
$ (b)$.
We show that the set of such $\delta$ is stationary in $S^\lambda_\theta$.

We now get to use 5. from 1.13.a.
Namely, 
$$E=\{\delta<\lambda:\,(\forall \gamma<\delta)\, (h(\gamma)<\delta)\}$$
is a club in $\lambda$. Therefore,
$$T=T_h\deq\{\delta\in S:\, \forall \alpha<\beta\,\,(
\alpha,\beta\in C_\delta\implies (\alpha,\beta]\cap E\neq\emptyset)\}$$
is stationary. It is easily seen that every element of $T$
satisfies $(b)$, so the
set of all $\delta\in S$ for which $j_\delta$ is defined,
is stationary.

Now, for some $i(\ast)<\kappa$, the set
$\{\delta\in S:j_\delta=i(\ast)\}$ is stationary.
Therefore,
$$\{\eta=\gamma^\ast(\delta,\epsilon):\,\epsilon\in S^\lambda_\kappa\,
\,\&\,\,
i_\delta(\epsilon)\hbox{ is defined and }=i(\ast)\}$$
is stationary in $\lambda$. For every such $\eta=\gamma^\ast(\delta,\epsilon)$,
we define $C=C(\eta)$ by $C\deq C_\delta\cap\gamma^\ast(\delta,\epsilon)$.

We can easily check that
this is a well posed definition and that
$i=i(\ast)$ is as required.\eop${}_{1.13.c.}$

\smallskip

{\it Continuation of the proof of 1.13.\/}
Now we apply 
the previous claim to $h_\xi\,(\xi<\lambda^+)$. For every $\xi
<\lambda^+$,
we fix $i(\xi)<\kappa$ as guaranteed by the claim. Then
note that for some $i(\ast)<\kappa$,
the set
$$W=\{\xi<\lambda^+: i(\xi)=i(\ast)\}$$
is unbounded in $\lambda^+$.
Of course, we can in fact assume $W$ to be stationary, but we
only need it to be unbounded.

We now define a new family of functions, based on the $h_\xi$'s.

For every $ \xi\in W$ and $\eta\in S^\lambda_\kappa$ let $h_{\xi,\eta}$
be the function with domain $a^\ast_\eta\deq\cup_{\alpha\in E_\eta}
a^\alpha_{i(\ast)}$, defined by
$$h_{\xi,\eta}(\beta)=\Min (a^\ast_\eta\setminus h_\xi(\beta)).$$
Observe that $a^{\ast}_\eta\subseteq\eta$
and $\eta=\sup(a^\ast_\eta)$. 
We have noted before that for $\zeta<\xi<\lambda$, we can fix an
$\alpha_{\zeta,\xi}\in\lambda$ such that
$$h_\zeta\rest[\alpha_{\zeta,\xi},\lambda)<h_\xi\rest[\alpha_{\zeta,\xi},\lambda).$$
So, if
$\zeta\le\xi\in W$, then
$$\alpha\in[\alpha_{\zeta,\xi},
\lambda)\implies h_{\zeta,\eta}(\alpha)\leq h_{\xi,\eta}(\alpha).$$
Also, if $\zeta<\xi\in W$, then
$$h_\zeta(\beta)<h_\xi(\beta)\wedge h_{\zeta,\eta}(\beta)=h_{\xi,\eta}(\beta)
\implies
[h_\zeta(\beta), h_\xi(\beta))\neq\emptyset\,\&\,
[h_\zeta(\beta), h_\xi(\beta))\cap a^\ast_\eta=\emptyset.$$
 
Using the above functions, we define the following sets $A_{\zeta,\xi}$ for
$\zeta<\xi\in W$.
$$A_{\zeta,\xi}\deq\{\eta\in S^\lambda_\kappa:\hbox{ for unboundedly
 } \beta\in a ^\ast_{\eta},
\hbox{ we have } h_{\zeta,\eta}(\beta) <h_{\xi,\eta}(\beta)\}.$$
We now show that these sets witness that $\DD_\lambda\rest S^\lambda_\kappa$
is not $\lambda^+$-saturated.

Note:

\item{(A)} If $\zeta<\xi\in W$, then $A_{\zeta,\xi}$ is a stationary
subset of $\lambda$.

[Why?
To see this, fix such $\zeta<\xi$ and suppose that $E$ is a club in $\lambda$
which misses $A_{\zeta,\xi}$. 

Suppose that $\eta>\alpha_{\zeta,\xi}$ and $\eta\in S^\lambda_\kappa\cap E$.
Then $\eta\notin A_{\zeta,\xi}$,
so there is $\beta_0<\sup (a^\ast_\eta)
=\eta$, such that for all $\beta\in
(a^\ast_\eta\setminus \beta_0)$, 
$$h_{\zeta,\eta}(\beta)=h_{\xi,\eta}(\beta),$$
or, equivalently
$$\bigl [h_\zeta(\beta), h_\xi(\beta)\bigr)\cap a^\ast_\eta=
\emptyset.$$
We can also assume $\beta_0\ge\alpha_{\zeta,\xi}$, so $\bigl
[h_\zeta(\beta),h_\xi(\beta)\bigl
)\neq\emptyset.$ In particular, $h_\zeta(\beta)\notin a^\ast_\eta$.
We can further 
assume that $\eta$ satisfies the conclusion of Claim 1.13.c,
with $h=h_\zeta$ and $i=i(\ast)$  (since $\zeta\in W$). Let $C$ be as there.

But then the conclusion of Claim 1.13.c tells us that
we can find a $\beta$ in $C$ which is greater than $\beta_0$ and
such that $h_\zeta(\beta)\in a^\ast_\eta$.]

\item{(B)} If $\zeta_1\leq\zeta_2\leq\xi_2\leq\xi_1\in W$, then $A_{\zeta_2,
\xi_2}\setminus A_{\zeta_1,\xi_1}$ is bounded.

[Why? This follows easily by the remarks after the definition of
$h_{\zeta,\eta}$.]

\smallskip

Now assume, for contradiction, that $\DD_\lambda\rest S^\lambda_\kappa$
is $\lambda^+$-saturated. Then, by (B):

\item{(C)} For each $\zeta\in W$, 
we have $\langle A_{\zeta,\xi}/\DD_\lambda:\,
\xi\in(\zeta,\lambda^+)\cap W
\rangle$ is eventually constant, say for
$\xi\in [\xi_\zeta, \lambda^+)$.

So, $A_\zeta\deq A_{\zeta,{\xi_\zeta}}$ satisfies $\zeta_1<\zeta_2
\in W\implies
A_{\zeta_1}\supseteq A_{\zeta_2}( \mod \DD_\lambda)$
(again by (B)). Hence, again by
the $\lambda^+$-saturation of $\DD_\lambda\rest S^\lambda_\kappa$,

\item{(D)} $\langle A_\zeta/\DD_\lambda:\,
\zeta\in W\rangle$ is eventually constant, say for $\zeta\geq\zeta^\ast$.
 
Choose $\zeta_\epsilon \in W\setminus\zeta^\ast$ for $\epsilon <\mu_{i(\ast)}^+$
such that
$$\epsilon(1)<\epsilon(2)\implies\xi_{\zeta_{\epsilon(1)}} <
\zeta_{\epsilon(2)},$$
which is possible since $W$ is unbounded.
By {\it (a)\/}
after the definition of $h_\xi$, we can find an $\alpha^\ast <\lambda$ such that
$$\epsilon(1)<\epsilon(2)<\mu_{i(\ast)}^+
\wedge \alpha^\ast\leq\alpha<\lambda\implies
h_{\zeta_{\epsilon(1)}}(\alpha)< h_{\xi_{\zeta_{\epsilon(1)}}}(\alpha)<
h_{\zeta_{\epsilon(2)}}(\alpha).\eqno{\oplus}$$
By clauses
(C) and (D), for all $\zeta,\xi \in W$ and $\epsilon<\mu^+_{i(\ast)}$,
$$A_{\zeta^\ast}=
A_{\zeta_\epsilon}=
A_{{\zeta_\epsilon},{\xi_{\zeta_\epsilon}}}=A_{{\zeta_\epsilon},
{\zeta_{\epsilon+1}}}\neq\emptyset\,(\mod\DD_\lambda).$$

By (C) and (D),
$$\cap_{\epsilon<{\mu_{i(\ast)}^+}} A_{{\zeta_\epsilon},
{\zeta_{\epsilon+1}}} = A_{\zeta^\ast}(\mod \DD_\lambda).$$
So,
$$A\deq\cap_{\epsilon<{\mu_{i(\ast)}^+}}
A_{{\zeta_\epsilon},{\zeta_{\epsilon+1}}}
\neq\emptyset\, (\mod \DD_\lambda).$$
Now we can choose an $\eta\in A
\setminus(\alpha^\ast+1),$ hence by $\oplus$,
$$\epsilon(1)<\epsilon(2)<\mu_{i(\ast)}^+\wedge\beta\in a_\eta^\ast\implies
h_{\eta,\epsilon(1)}(\beta)\leq h_{\eta,{\zeta_{\epsilon(2)}}}(\beta).
\eqno{\oplus_1}$$
For each $\beta\in a _{\eta^\ast}\setminus(\alpha^\ast+1)$, the sequence
$\langle h_{\eta,\zeta_\epsilon}(\beta):\,\epsilon<\mu^+_{i(\ast)}\rangle$
is non-decreasing. Hence,
since $\card{a^\ast_\eta}\le\mu_{i(\ast)}$,
for some $\epsilon(\beta)<\mu^+_{i(\ast)}$,
the sequence
$\langle h_{\eta,\zeta_\epsilon}(\beta):\,\epsilon(\beta)\leq
\epsilon<\mu^+_{i(\ast)}\rangle$ is constant. Let $\epsilon(\ast)\deq
\sup_{\beta\in a_\eta^\ast\setminus(\alpha^\ast+1)}
(\epsilon(\beta))<\mu^+_{i(\ast)}$.
But $\eta\in A_{\zeta_{\epsilon(\ast)},\zeta_{\epsilon(\ast)+1}}$,
hence by the definition of $A_{\zeta,\xi}$'s,
there is a $\beta\in 
a^\ast_\eta\setminus(\alpha^\ast+1)$,
such that
$$h_{\zeta_{\epsilon(\ast),\eta}}(\beta)<h_{\zeta_{\epsilon(\ast) +1},
\eta} (\beta).
\eqno{\oplus_2}$$
So we get a contradiction.
\eop${}_{1.13.}$ 

\medskip

A similar argument can be applied to other normal filters
$\DD$ on
$\lambda$, under certain conditions. 
In addition we shall see
that, under some assumptions on the cardinal arithmetic,
the fact that $\DD$ is not $\lambda^+$-saturated
is strongly witnessed, by the existence of a $\diamondsuit$
on $\DD$. That is, we obtain $\diamondsuit_{\DD}^\ast
(S^\lambda_\kappa)$. This notation is explained
in the following

\smallskip

{\bf Definition 1.14.} If $\DD$ is a normal filter on $\lambda$,
and $S$ is $\DD$-stationary, then $\diamondsuit_{\DD}(S)$ means:

There is a sequence $\langle A_\alpha:\,\alpha\in S\rangle$
such that each $A_\alpha\subseteq\alpha$ and for every $A\in
[\lambda]^\lambda$,
$$\{\alpha\in S:\,A\cap\alpha=A_\alpha\}\in \DD^+.$$
The statement $\diamondsuit_{\DD}^\ast(S)$
means:

There is a sequence $\langle\PP_\alpha:\,\alpha\in S\rangle$
such that each $\PP_\alpha\subseteq\PP(\alpha)$ and $\card{\PP_\alpha}
\leq\alpha$, and for every $A\subseteq\lambda$, there is a $C\in \DD$ such
that for all $\alpha$,
$$\alpha\in S\cap C\implies A\cap\alpha\in\PP_\alpha.$$

\smallskip

By a well known result of Kunen (see [Ku 2]),
$$\diamondsuit_{\DD}^\ast(S)\implies\diamondsuit_{\DD}(S).$$
It is easily seen that $\diamondsuit_{\DD}(S)$ implies the existence of
an almost disjoint family of $\DD$-stationary subsets of $\lambda$, of
size $2^\lambda$. Therefore, if $\diamondsuit_{\DD}(S)$
holds, $\DD$ is not $2^\lambda$-saturated.

Looking back at the proof of Theorem 1.13,
there are two important facts that we were using. The first is that there
is a $\theta\in(\kappa,\lambda)$ such that $\theta$ has the weak
reflection property for $\kappa$. 
The other important ingredient of the proof
is Fact 1.13.c.
Dealing with filters other than $\DD_\lambda$,
to obtain the corresponding
version of 1.13.c, we have to strengthen our
assumptions,
Here, $\II[\theta,\kappa)$ is as in the Definition 1.3.

\medskip

{\bf Theorem 1.15.} Assume $\mu>\theta=\cf(\theta)>\kappa=\cf(\mu)>\aleph_0$,
while
$\lambda=\mu^+$ and $S^\theta_\kappa\notin\II [\theta,\kappa)$. Suppose
that $\DD$
is a normal filter on $\lambda$ such that $S^\lambda_\kappa\in\DD$.

If, for some $S$,

\item{$(\ast)$} For every $A\in\DD$,
$$\Bigl\{\delta\in S: \forall\, C\hbox{ a club in }\delta\,\bigl
(\{\otp(C\cap\gamma):\,
\gamma\in C\cap A\}\neq\emptyset\mod \II[\theta,\kappa)\bigr)\Bigr\}$$
is a $\DD$-stationary set in $\lambda$,

{\it then:\/}

\item{(1)} Claim 1.13.c holds for
any $\bar{C}$ and $a^\delta_i$ as in the assumptions of the Claim 1.13.c,
with ``stationary $\subseteq\lambda$''
replaced by ``$\DD$-stationary $\subseteq\lambda$''.
\item{(2)} $\DD\rest S^\lambda_\kappa$ is not $\lambda^+$-saturated.
\item{(3)} If $2^\mu=\lambda$ and $\mu^{[\kappa]}=\mu$, then
$\diamondsuit^\ast_{\DD}(S^\lambda_\kappa)$ holds.

\medskip

We remind the reader of the notation
$\mu^{[\theta]}$ for the revised cardinal power, 
from [Sh 460].
 
\smallskip

{\bf Definition 1.16.} Suppose that $\nu>\chi$ are infinite cardinals
and $\chi$ is regular. Then
$$\nu^{[\chi]}=\Min\{\card{\PP}:\,\PP\subseteq [\nu]^\chi\hbox{ and }
\forall A\in [\nu]^\chi\,
(A\subseteq\hbox{ the union of} <\chi\hbox{ elements of }\PP)\}.$$

\smallskip

{\bf Proof of 1.15.} (1) and (2) are
easily adjusted from Theorem 1.13.

(3) The conclusion is the same as that of 
[Sh 186,\S3], and the proof is the same.
The assumptions on the
cardinal arithmetic are here the same as in [Sh 186,\S3],
with $\chi$ from there equal to our $\kappa$.
The only difference is that
the proof in [Sh 186,\S3] started from
$\square\langle\lambda\rangle$, but we can use
Fact 1.13.a instead.\eop${}_{1.15.}$

\bigskip

{\bf \S2. The $\club^\ast_{-\bar\theta}(S)$
principle.}
Suppose that $S$ is a stationary
subset of a regular uncountable cardinal $\lambda$.
The aim of this section is to formulate
a combinatorial principle which 
suffices to show that certain normal filters on
$\lambda$ cannot be $\lambda^+$-saturated.
The principle,
$\clubsuit^\ast_{-{\bar\theta}}(S)$,
is a form of $\clubsuit(S)$, where 
${\bar\theta}$ is a sequence of ordinals.
Our interest in this comes from two facts
presented in
2.8.
Firstly,
if $\lambda$ is a successor of a singular cardinal
$\mu$
of uncountable cofinality
$\kappa$, then $\clubsuit^\ast_{-\mu}(
\lambda\setminus S^\lambda_\kappa)
$ is always true. This can be used
to obtain an alternative proof of that part of the result from
[Sh 212, 14]=[Sh 247, 6] 
which states that no normal filter concentrating on $\lambda\setminus
S^\lambda_\kappa$ is $\lambda^+$-saturated.
Secondly, if $\II[\lambda,\kappa)\rest S^\lambda_\kappa$ contains
only nonstationary sets, $\club^\ast_{-\mu}(S^\lambda_\kappa)$ is true,
so by 2.5. we can conclude that $\DD\rest S^\lambda_\kappa$ is
not $\lambda^+$-saturated.
The key to the proof of 2.8. is the combinatorial lemma 2.7. 

We commence by recalling the definition of some versions of $\club$.

\smallskip

{\bf Definition 2.0.} Let
$\lambda$ be an uncountable regular cardinal
and $S\subseteq\lambda$ stationary. Then:

\smallskip

(0) $\clubsuit(S)$ means:

There is a sequence $\bar{A}=\langle A_\alpha:\,\alpha\in S\cap LIM\rangle$
such that:

\item{{\it (i)\/}} For each $\alpha\in S\cap LIM$, the set $A_\alpha$
is an unbounded subset of $\alpha$.
\item{{\it (ii)\/}} For
all
$A\in[\lambda]^\lambda$, the
set
$$G_{\clubsuit(S)}^{\bar{A}}[A]\deq\{
\alpha\in S \cap LIM:\,
A_\alpha\subseteq A\}$$
is non-empty.

\smallskip

(1) $\clubsuit^\ast (S)$ means:

There is a sequence $\bar\PP=\langle\PP_\alpha:\,\alpha
\in S\cap LIM\rangle$ such that

\item{{\it (i)\/}} For each $\alpha\in S\cap LIM$, we have
$\card{\PP_\alpha}\le\card{\alpha}$.
\item{{\it (ii)\/}} If $B\in \PP_\alpha$, then $B$ is an unbounded
subset of $\alpha$.
\item{{\it (iii)\/}} For every $A\in[\lambda]^\lambda$,
there is a club $C_A$ of $\lambda$
such that
$$G_{\club^\ast(S)}^{\bar{\PP}}[A]\deq\{\alpha
\in S\cap LIM:\,\exists B\in\PP_\alpha\,(B\subseteq A)\}\supseteq
C_A\cap S\cap LIM.$$

\medskip

Many authors use a different definition of $\club$, in which every
unbounded set is required to be ``guessed'' stationarily many times.
The following well known
fact shows that the two definitions are equivalent.
We also include some other
easy observations about Definition 2.0.

\medskip

{\bf Fact 2.1.} Assume that $\lambda$ and $S$ are as in Definition 2.0.

(0) If
$\bar{A}=\langle A_\alpha:\,\alpha\in S\cap LIM
\rangle$ is a $\club(S)$
-sequence,
{\it then\/},
for every $A\in[\lambda]^\lambda$, the set $G_{\clubsuit(S)}^{\bar{A}}[A]$
is stationary.


(1) If $\club(S)$ 
holds,
{\it then\/} there is a $\club(S)$-sequence
$\langle A_\alpha:\,\alpha\in S\cap LIM\rangle$ such that for each $\alpha\in
S
\cap LIM$, we have that $\otp(A_\alpha)=\cf(\alpha).$



(3) Suppose that $\DD$ is a normal filter on $\lambda$ and 
$S\in\DD^+$ is such that 
$\bar{\PP}$ exemplifies that $\clubsuit^\ast(S)$ holds. {\it Then,\/}
for every $A\in[\lambda]^\lambda$, the set $G_{\clubsuit^\ast(S)}^{\bar{\PP}}
[A]$ is $\DD$-stationary.


{\bf Proof.} (0) Otherwise, we could find an $A\in[\lambda]^\lambda$
and a club $C$ in $\lambda$, such that for all $\alpha$ in $C\cap S\cap LIM$,
the set $A_\alpha$
is not a subset of $A$.
Then set $A^{\dagger}
=\{\Min(A\setminus\alpha):\alpha\in C\}$, so
$A^\dagger\in[\lambda]^\lambda$.
But if $A_\alpha\subseteq A^\dagger$, then,
$A_\alpha$ is
also a subset of $A$.
On the other hand,
since $A^\dagger$ is unbounded in $\alpha$ (as $A_\alpha$ is),
also $C$ is unbounded in $\alpha$, therefore $\alpha\in
C$. This is a contradiction.

(1) Suppose that $\langle B_\alpha:\,
\alpha\in S\cap LIM\rangle$ exemplifies
$\club(S)$ 
and define
for each $\alpha\in S\cap LIM$, the set $A_\alpha$ to be any cofinal
subset of $B_\alpha$ with $\otp(A_\alpha)=\cf(\alpha)$.

(3) We simply remind the reader of the following elementary

{\it Observation.\/} For every 
club $C$ of $\lambda$, we have $C\in\DD$ (so the set
$S\cap C$ is $\DD$-stationary).

[Why? Suppose that $C$ is a club of $\lambda$ such that
$C\notin\DD$, so
$S\setminus C\in\DD^+$.
We define the following function, for
$\alpha\in S\setminus C$:
$$f(\alpha)\deq\sup(C\cap\alpha)$$
and we note that $f$ is regressive on $S\setminus C$.
Then we can find a $T\subseteq
S\setminus C$
which is $\DD$-stationary
and such that $f\rest T$ is a constant.
Then $T$ must be bounded, which is a contradiction.]

\eop${}_{2.1.}$

\medskip

Now
we introduce the
version of the $\club^\ast$ principle that will mainly interest
us.
The guessing requirement is weaker, while
the order type of the sets entering each family in the $\club^\ast$-sequence
is controlled by a sequence of ordinals.

\medskip

{\bf Definition 2.2.} Let $\lambda$
be a regular uncountable cardinal
and $S\subseteq
\lambda$ stationary.
Let, $\bar\theta=\langle\theta_\alpha:\,\alpha\in S\rangle$ be
a sequence of ordinals. Then

\smallskip

$\clubsuit^\ast_{-\bar\theta}(S)$
means:

There is a sequence $\bar{\PP}=\langle \PP_\alpha:\,\alpha\in S\cap LIM\rangle$
such that:

\item{{\it (i)}}
For each $\alpha\in S\cap LIM$, 
the family
$\PP_\alpha$ consists of $\le\card{\alpha}$ unbounded subsets
of $\alpha$.

\item{{\it (ii)}} For each $\alpha\in S\cap LIM$
and $B\in\PP_\alpha$, we have
$\otp(B)<\theta_\alpha$.

\item{{\it (iii)}} For $A\in[\lambda]^\lambda$,
there is a club $C_A$ of $\lambda$ such that
$$G_{\clubsuit^\ast_{-\bar{\theta}}(S)}^{\bar{\PP}}[A]\deq\{\alpha:\,
\bigl(\exists B\in\PP_\alpha\bigr)\bigl(\alpha=
\sup(B\cap A)\bigr)\}\supseteq C_A\cap S\cap LIM.$$

\smallskip

If for each $\alpha\in S$, $\theta_\alpha=
\cf(\alpha)+1$,
we omit
$\bar\theta$ in the above notation.

If for some $\mu$ we have that
$\theta_\alpha=\mu$ for all $\alpha\in S\cap LIM$, then
we write $\club^\ast_{-\mu}(S)$ rather than $\club^\ast_{-\bar\theta}(S)$.

\medskip

We make some easy remarks on Definition 2.2.

\medskip

{\bf Observation 2.3.} Assume that $\lambda$, $S$
and $\bar\theta$ are as in Definition
2.2. 

(0) If the set of all $\alpha\in S\cap LIM$ for which $\theta_\alpha>
\cf(\alpha)$ is non-stationary, {\it then\/}
$\clubsuit^\ast_{-\bar\theta}(S)$
is false. Otherwise $\club^\ast(S)\implies \club^\ast_{-\bar\theta}(S)$.

(1) If
$\Sigma_{\gamma<\theta_\alpha}\card{\gamma}^{\cf(\alpha)}\leq\card{\alpha}$,
for each $\alpha\in S$, {\it then\/}
$$\clubsuit^\ast_{-\bar\theta}(S)\implies \clubsuit^\ast(S).$$

(2) Suppose that $\DD$ is a normal filter on $\lambda$, while $S$ is
$\DD$-stationary and $\bar{\PP}$ exemplifies $\club^\ast_{-\bar\theta}(S).$
{\it Then\/}, for every $A\in[\lambda]^\lambda$
the set $G_{\club^\ast_{-\bar\theta}(S)}^{\bar\PP}[A]$ is $\DD$-stationary.




{\bf Proof.} (0) Obvious.

(1) If $\bar\PP=
\langle\PP_\alpha :\,\alpha\in S\cap LIM\rangle$ exemplifies $\clubsuit^\ast_{
-\bar\theta}(S)$,
define 
$$\PP_\alpha^{\dagger}=
\cases{
\Bigl\{ B:\exists A\in\PP_\alpha\,\bigl(B\subseteq A
\,\wedge B\hbox{ cofinal in }
\alpha\,\wedge\, \otp(B)=\cf(\alpha)\bigr)\Bigr \} & if
$\theta_\alpha>\cf(\alpha)$
\cr
\{\alpha\} & otherwise.\cr}
$$
Then 
$\card{\PP_\alpha^\dagger}=\card{\PP_\alpha}^{\cf(\alpha)}\le\alpha$
and, by (0),
$\langle\PP_\alpha^\dagger:\,\alpha\in S\cap LIM\rangle$
exemplifies $\clubsuit^\ast(S)$.

(2) Like 2.1.3.
\eop${}_{2.3.}$

\medskip

We can consider also $\club$-sequences whose failure to guess
is always confined to a set in a given ideal on $\lambda$.
In this context, for example $\clubsuit(S)$ will mean
$\clubsuit(S)/J_\lambda$.

\medskip

{\bf Definition 2.4.} Let $\lambda$ and $S$ be as above, while
$\II$ is an ideal on $\lambda$. Then

$\club(S)/\II$ means:

There is a sequence $\bar{A}=
\langle A_\alpha:\,\alpha\in S\rangle$
such that each $A_\alpha$ is 
cofinal
subset of $\alpha$,
and this sequence has the
following property. For every $A\in[\lambda]^\lambda$,
the set
$$G_{\club(S)/\II}^{\bar{A}}
[A]\deq\{\alpha:
\,A_\alpha\subseteq A\}$$
satisfies $G_{\club(S)/\II}^{\bar{A}}[A]\notin \II$.

If $\DD$ is the dual filter of the ideal $\II$, then $\club(S)/\DD$ means the
same as $\club(S)/\II$.
We extend this definition in the obvious way to the other mentioned versions
of the $\club$ principle.

\medskip

{\bf Theorem 2.5.} Suppose that $\lambda=\mu^+$ and $\mu$ is a limit
cardinal.

(1) Assume that $\club^\ast_{-\mu}(S)$ holds for a stationary $S\subseteq
\lambda$. {\it Then\/} no normal filter $\DD$ on $S$ is $\lambda^+$-saturated.

(2) If $S\subseteq\lambda$ is stationary and $\DD$ is a normal filter on $S$
such that $\club^\ast_{-\mu}(S)/\DD$ holds, {\it then\/} $\DD$ is not
$\lambda^+$-saturated.

{\bf Proof.} Let us fix a sequence $\bar\PP=\langle
\PP_\delta:\,\delta\in S\cap LIM\rangle$ which exemplifies $\club^\ast_{-
\mu}(S)$. Therefore, for each $\delta\in S\cap LIM$, we have a family
$\PP_\delta=\{P_{\delta,i}:\,i<i^\delta\le\delta\}$ such that
each $P_{\delta,i}$ is a cofinal subset of $\delta$ of order type
$\otp(P_{\delta,i})<\mu$. 

\def\pr{\rm{pr}}

We fix a
1-1 onto pairing function $\pr:\lambda\times\lambda\to\lambda
\setminus\omega$ such that
for each $\alpha,\beta\in\lambda$, we have
$$\Max\{\alpha,\beta\}\le \pr(\alpha,\beta)<(\card{\alpha}+\card{\beta})^+.$$
(Here, we use the convention that $n^+=\aleph_1$ for $n\in\omega$.)

We shall also fix a club $C$ of $\lambda\setminus\mu$ such that
$$\alpha,\beta<\gamma\,\,\&\,\,\gamma\in C\implies \pr(\alpha,\beta)<\gamma.$$
Now, we choose sets $b^\zeta_\alpha$ for $\zeta<\lambda^+$ and $\alpha
<\lambda$, and the functions $h_\zeta:\lambda\to\lambda$ for $\zeta<
\lambda^+$ as in the proof of 1.13, and with the same properties as there.
For $\zeta<\lambda^+$, we define the unbounded subset $X_\zeta$ of
$\lambda$ by
$$X_\zeta\deq\bigl\{\hbox{pr}\bigl(\alpha,h_\zeta(\alpha)\bigr):\,
\alpha<\lambda\bigr\},$$
and partial functions $g_\zeta$ by
$$g_\zeta(\delta)\deq\Min\{i<\delta:\,\sup(X_\zeta\cap P_{\delta,i})=
\delta\}.$$
In fact, the domain of each $g_\zeta$ is exactly the set
$G_{\club^\ast_{-\mu}(S)}^{\bar\PP}[X_\zeta]$, so $\hbox{Dom}(g_\zeta)$
is $\DD$-stationary (by 2.3.2). For $\delta\in {\Dom}
(g_\zeta)\cap C$, we can
define
$$f_\zeta(
\delta)\deq\hbox{pr}\bigl( g_\zeta(\delta),{\card{P_{\delta,g_\zeta(\delta)}}}^+
\bigr).$$
Notice that $g_\zeta$
is regressive on its domain, so if $\delta\ge\mu
$ and $f_\zeta(\delta)$ is defined, then $f_\zeta(\delta)<\delta$
(as $\delta\in C$).
Therefore, there is a $\DD$-stationary
set $B_\zeta$ such that $f_\zeta\rest B_\zeta$ is constantly equal
to $\hbox{pr}(i_\zeta,\theta_\zeta)$ for some $i_\zeta$ and a regular
cardinal $\theta_\zeta<\mu$. Then there are $i^\ast<\lambda$ and a
regular cardinal $\theta^\ast<\mu$ such that the set
$$W\deq
\bigl\{\zeta<\lambda^+:\,\bigl(i_\zeta,\theta_\zeta
\bigr)=\bigl(i^\ast,\theta^\ast
\bigr)\bigr\}$$
is unbounded.

Let us define for $S\cap \acc(C)\cap LIM$
with $i^\ast<i^\delta$ the set
$$d_\delta\deq\bigl\{\alpha:\,\bigl(\exists\beta\bigr)\bigl(\hbox{pr}
(\alpha,\beta)\in P_{\delta,i^\ast}\bigr)\bigr\}\cup
\bigl\{\beta:\,\bigl(\exists\alpha\bigr)\bigl(\hbox{pr}(\alpha,\beta)
\in P_{\delta,i^\ast}\bigr)\bigr\}.$$
By the definition of $\pr$, we have $d_\delta\subseteq\delta$. In fact,
since $\delta\in\acc(C)$, the set
$d_\delta$ is an unbounded subset of $\delta$. Like in 1.13, we can
define for $\zeta\in W$ and $\delta$
for which $d_\delta$ is defined, a function
$h_{\zeta,\delta}$ on $d_\delta$ given by
$$h_{\zeta,\delta}(\beta)=\cases{\Min\bigl(d_\delta\setminus
(h_\zeta(\beta)+1)\bigr) & if $d_\delta\setminus \bigl(h_\zeta(\beta)+1\bigr)
\neq\emptyset$\cr
\lambda & otherwise.\cr}$$
Therefore
$$h_{\zeta,\delta}:\,d_\delta\to d_\delta\cup\{\lambda\}.$$
For $\zeta,\xi\in W$ we define sets
$$A_{\zeta,\xi}=\{\delta\in S:\,d_\delta\hbox{
is defined and for unboundedly }
\beta\in d_\delta \hbox{ we have }h_{\zeta,\delta}(\beta)
<h_{\xi,\delta}(\beta)\}.$$
Assume now that $\DD$ is $\lambda^+$-saturated, and let us make some simple
observations about the just defined sets:

\item{\it{(a)\/}} $A_{\zeta,\xi}/\DD$ increase with $\xi$ and decrease
with $\zeta$.
\item{\it{(b)\/}} Since $\DD$ is $\lambda^+$-saturated, for any fixed $\zeta
\in W$, the sequence 
$\langle A_{\zeta,\xi}/\DD:\,\xi<\lambda^+\rangle$
is eventually constant, let us say for $\xi\in[\xi_\zeta,\lambda^+)
\cap W$.

Similarly,

\item{\it{(c)\/}} $\langle A_\zeta\deq A_{\zeta,\xi_\zeta}/\DD:\,
\zeta<\lambda^+\rangle$ are eventually constant, say for $\zeta\in W\setminus
\zeta(\ast).$

We choose by induction on $\epsilon <\theta^\ast$ ordinals
$\zeta_\epsilon$ such that

\item{\it{(i)\/}} $\zeta_\epsilon\in\bigl(\zeta(\ast),\lambda^+\bigr).$
\item{\it{(ii)\/}} $\zeta_\epsilon\in W$.
\item{\it{(iii)\/}} $\zeta_\epsilon$ are strictly increasing with $\epsilon$.
\item{\it{(iv)\/}} $\zeta_{\epsilon+1}>\xi_{\zeta_\epsilon}.$

Therefore
$$\epsilon(1)<\epsilon(2)<\theta^\ast\implies A_{\zeta_{\epsilon(1)},
\zeta_{\epsilon(2)}}/\DD=A_{\zeta_{\epsilon(1)},\xi_{\zeta_{\epsilon(1)}}}/\DD
=A_{\zeta(\ast)}/\DD.$$
Let $\gamma(\ast)<\lambda$ be such that for all $\alpha\in \bigl(\gamma
(\ast),\lambda^+\bigr)$, the sequence $\langle h_{\zeta_\epsilon}
(\alpha):\,\epsilon<\theta^\ast\rangle$ is strictly increasing. (Recall
that for $\zeta<\xi<\lambda^+$, we know that $h_\zeta$ is eventually
strictly less than $h_\xi$.) We can as well assume that $\gamma(\ast)>\omega$.

By the above and the fact that $\theta^\ast<\lambda$, we can find
a set $E\in\DD$ such that
$$\epsilon(1)<\epsilon(2)<\theta^\ast\implies A_{\zeta_{\epsilon(1)},
\zeta_{\epsilon(2)}}\cap E= A_{\zeta_{\epsilon(1)},\xi_{\zeta_{\epsilon(1)}}}
\cap E= A_{\zeta(\ast)}\cap E.$$
We can assume that $\Min(E)>\gamma(\ast).$ Without loss of generality,
we can also add that $E\subseteq \acc(C)$ and
$$\epsilon <\theta^\ast\,\&\,\delta\in E\cap S\cap LIM\implies
h_{\zeta_\epsilon}{}^{''}\delta\subseteq\delta,$$
since $\theta^\ast<\lambda$.

Now we discuss the two possible cases, the first of which corresponds
to the situation in 1.13.

\smallskip

$\underline{\hbox{Case 1.}}$ 
For some $\bar\zeta\ge\zeta(\ast)$,
we have $\bar\zeta\in W$ and $A_{\zeta(\ast)}\cap B_{\bar\zeta}$
is $\DD$-stationary.

We choose a $\delta\in E\cap
B_{\bar\zeta}\cap A_{\zeta(\ast)}$.
In particular, $d_\delta$ is defined. We consider
$\langle h_{{\zeta_\epsilon},\delta}:\,\epsilon<\theta^\ast\rangle$.
By the choice of $E$, each $h_{{\zeta_\epsilon},\delta}$ is a function
from $d_\delta$ to itself.

For any $\alpha\in d_\delta\setminus \gamma(\ast)$, the sequence
$\langle h_{{\zeta_\epsilon},\delta}(\alpha):\,\epsilon<\theta^\ast\rangle$
is non-decreasing. As
$\delta\in B_{\bar\zeta}$ and
$\bar\zeta\in W$, we know that $\card{d_\delta}\leq\card{P_{\delta,i^\ast}}<
\theta^\ast$, so
the above sequence is eventually constant.

Similarly, $\card{d_\delta\setminus\gamma(\ast)}<\theta^\ast$,
so there is some $\epsilon_0<\theta^\ast$ and some function $h$
such that for all $\epsilon\in(\epsilon_0,\theta^\ast)$, we have
$h_{\zeta_\epsilon,\delta}\rest\bigl(d_\delta
\setminus\gamma(\ast)\bigr)=h.$
But
$\delta\in A_{\zeta(\ast)}\cap E$, so $\delta\in A_{{\zeta_\epsilon},
{\zeta_{\epsilon+1}}}$ for all $\epsilon<\theta^\ast$, and therefore
for any such $\epsilon$,
$$\{\alpha\in d_\delta:\,h_{{\zeta_\epsilon},\delta}(\alpha)
<h_{{\zeta_{\epsilon+1}},\delta}(\alpha)\}$$
is unbounded in $\delta$. This is a contradiction.

\smallskip

$\underline{\hbox{Case 2.}}$ $A_{\zeta(\ast)}
\cap B_{\bar\zeta}$ is not $\DD$-stationary for any $\bar\zeta\in W$
with $\bar\zeta\ge\zeta(\ast)$.

Then for all $\zeta\in W$
and $\bar\zeta\in W\setminus \zeta(\ast)$, the set $A_{\zeta,{\xi_\zeta}}
\cap B_{\bar\zeta}$ is
not $\DD$-stationary, as $A_{\zeta,{\xi_\zeta}}\subseteq A_{\zeta(\ast)}
/\DD$. Similarly, since for $\zeta<\xi\in W$, we have that
$A_{\zeta,\xi}\subseteq A_{\zeta,{\xi_\zeta}}/\DD$, we conclude that
$A_{\zeta,\xi}\cap B_{\bar\zeta}$ is not $\DD$-stationary
for any $\zeta<\xi\in W$ and $\bar\zeta\in W\setminus\zeta(\ast)$.

On the other hand, as each $B_\zeta$ for $\zeta\in W$ is $\DD$-stationary,
and we are assuming that $\DD$ is $\lambda^+$-saturated, there are
$\zeta<\xi\in W
\setminus \zeta(\ast)$ such that $B_\zeta\cap B_\xi$ is $\DD$-stationary.
We fix such $\zeta$ and $\xi$.

Let us choose a $\delta\in (E\cap B_\zeta\cap B_\xi)\setminus A_{\zeta,\xi}$.
Without loss of generality, we can assume that there is an $\alpha_{\zeta,\xi}
<\delta$ such that $\alpha\ge\alpha_{\zeta,\xi}\implies h_{\zeta}(\alpha)
<h_\xi(\alpha).$
Note that $d_\delta$ is defined. Then, by the definition of
$A_{\zeta,\xi}$, we can find a $\gamma_1
\in\bigl(\alpha_{\zeta,\xi},\delta\bigr)$
such that
$$h_{\zeta,\delta}\rest(d_\delta\setminus\gamma_1)=
h_{\xi,\delta}\rest(d_\delta\setminus\gamma_1).$$
Note now that there is a club $E_\xi$ such that
$$\alpha\in E_\xi\implies h_\xi{}^{''}\alpha\subseteq\alpha,$$
and a similarly defined club $E_\zeta$. We can without loss of
generality assume that $\delta\in E_\zeta\cap E_\xi$, so both
$h_\zeta$ and $h_\xi$ are functions from $d_\delta$ to itself.
Now, since $\delta\in B_\zeta\cap B_\xi$, we know that $g_\zeta(\delta)=
g_\xi(\delta)=i(\ast)$. By the definition of $X_\xi$, we can find an
$\alpha\in(\gamma_1,\delta)$ such that ${\rm pr}\bigl(\alpha_1,
h_\xi(\alpha_1)\bigr)\in d_\delta$, so by the definition of $h_{\xi,\delta}$
we have
$$h_{\xi,\delta}(\alpha_1)>h_\xi(\alpha_1).$$
On the other hand, $h_\zeta(\alpha_1)<h_{\xi}(\alpha_1)\in d_\delta$, so
$h_{\zeta,\delta}(\alpha_1)\le h_\zeta(\alpha_1)$. Therefore,
$h_{\zeta,\delta}(\alpha_1)\neq h_{\xi,\delta}(\alpha_1)$, which is
a contradiction to $\alpha_1>\gamma_1$.

(2) Follows from the proof of (1).\eop${}_{2.5.}$

\medskip

Like $\diamondsuit$ and by the same proof, the usual $\clubsuit$
principle on $\lambda$ can be used for guessing not just unbounded subsets of
$\lambda$, but any other structure which can be coded by the unbounded
subsets of $\lambda$. With the $\clubsuit^\ast_-$ principle, this does
not seem to be the case, or at least the $\diamondsuit$-like proof fails.
In particular, we do not know if the following is true:

\medskip

{\bf Question 2.6.} Suppose that $\lambda$ is a regular uncountable
cardinal and $S$ a stationary subset of $\lambda$ such that
$\club^\ast_-(S)$ holds. Is it true that there is a sequence
$\langle \FF_\alpha:\,\alpha\in S\cap LIM\rangle$ with the
following properties:

\item{\it{ (i)}\/} Each $\FF_\alpha$ consists of partial functions from
$\alpha$ to $\alpha$, each of which has an unbounded subset of $\alpha$
as its domain.
\item{\it{ (ii)}\/} $\card{\FF_\alpha}\le\alpha.$
\item{\it{ (iii)}\/} For every function $f:\,\lambda\to\lambda$, there
is a club $C_f$ of $\lambda$ with the property
$$\alpha\in C_f\cap S\cap LIM\implies\bigl(\exists g\in\FF_\alpha\bigr)
\bigl(\sup\{\beta\in{\rm Dom}(g):\,g(\beta)=f(\beta)\}=\alpha\bigr)?$$

\medskip

We note that a positive answer to 2.6. would quite simplify the proof
of 2.5.

Our next goal is to prove that $\club^\ast_{-\mu}
(\lambda\setminus S^\lambda_\kappa)$ is true for any $\lambda$
which
is the successor of a singular cardinal $\mu$ of uncountable
cofinality $\kappa$. The key is the following

\medskip

{\bf Lemma 2.7.}
Assume that $\lambda=\mu^+$ and $\aleph_0<\kappa=\cf(\mu)
<\mu$. Also $\mu=\Sigma_{i<\kappa}\mu_i$, where $\langle\mu_i:\,i<\kappa
\rangle$ is an increasing continuous sequence of cardinals, and for simplicity,
$\mu_0>\kappa$.

 {\it Then\/}
there is a sequence
$$\big\langle \langle a^\alpha_i:\,i<\kappa,
\alpha<\lambda\bigr\rangle$$
such that for every $\alpha<\lambda$, the
sets $a_i^\alpha (i<\kappa)$ are subsets of $\alpha$
which are
$\subseteq$-increasing in $i$, with $\card{a_i^\alpha}\leq\mu_i$, and such that:

For any $f:\lambda\to\lambda$, if
$$A^f_i=\Bigl\{\alpha<\lambda:\,\alpha=\sup\{\zeta\in a_i^\alpha:\,f(\zeta)\in
a_i^\alpha\}\Bigr\}$$
{\it then\/}

\item{(A)} $A^f_i$ are $\subseteq$-increasing in $i$.
\item{(B)} $\lambda\setminus S^\lambda_{\kappa}\subseteq\cup_{i<\kappa}A^f_i
(\mod\DD_\lambda).$
\item{(C)} If $\gamma\in S_\kappa^\lambda$,
while $ i<\kappa$ and $A^f_i$ reflects on $\gamma$,
then $\gamma\in A^f_i$.
(In fact, this is true for any $\gamma<\lambda$ with $
\kappa\ge\cf(\gamma)>\aleph_0$.)

So

\item{(C${}^{'}$)} If $S_f\deq S^\lambda_\kappa\setminus\cup_{i<\kappa}
A^f_i$, then for no $i<\kappa$ does $A^f_i$ reflect in any $\delta
\in S_f$.
\item{(D)} There is a nonstationary set $N$
such that $S_f\setminus N\in \II[\lambda,\kappa).$

{\bf Proof.} We describe the choice of sets $a^\alpha_i$,
and then we check that all claims of the lemma are satisfied. 

First we fix a sequence $\langle e_\gamma:\,\gamma\in LIM\cap\lambda\rangle$
such that $e_\gamma$ is a club of $\gamma$ with $\otp(e_\gamma)=\cf(\gamma).$
Then we define by induction on $\alpha<\lambda$ sets $a^\alpha_i$ for all
$i<\kappa$, requiring that for all limit $\gamma$
$$\mu_i\ge\cf(\gamma)\implies \cup_{\beta\in e_\gamma}a_i^\beta\subseteq
a_i^\gamma,$$
in addition to the requirements that we already mentioned in the statement
of the lemma.

Now we suppose that $f:\lambda\to\lambda$ is given, and
check claims (A)--(D).

(A) This follows from the fact that for every $\alpha<\lambda$,
$$i<j<\kappa\implies a^\alpha_i\subseteq a^\alpha_j.$$

(B) Let
$$E\deq\Bigl\{\alpha\in\lambda:\,\alpha=\sup\{\zeta<\alpha:\,f(\zeta)<\alpha\}
\Bigr\}.$$
Then $E$ is a club of $\lambda$. We shall show that
$$(\lambda\setminus S^\lambda_\kappa)\cap\acc(E)\subseteq \cup_{i<\kappa}
A^f_i.$$
So, let us take a $\gamma\in\acc(E)$ such that $\cf(\gamma)\neq\kappa$.
Then $e_\gamma\cap E$ is a club of $\gamma$, and $\otp(
e_\gamma\cap E)=\cf(\gamma)$.
Let $e_\gamma
\cap E=\{\beta_\epsilon:\,\epsilon<\cf(\gamma)\}$ be an increasing
enumeration. So, for all $\epsilon$, we have
$$f(\beta_\epsilon)<\beta_{\epsilon+1}<\gamma.$$
Since $\cf(\gamma)\neq\kappa$, there must be an $i<\kappa$
such that
for some unbounded $c\subseteq
\cf(\gamma)$, we have $\{\beta_\epsilon, f(\beta_\epsilon):\,\epsilon\in c\}
\subseteq a_i^\gamma$ (note that we are using the fact that $a^\alpha_i$
are increasing with $i$).
Then $\gamma\in A^f_i$.

(C) Suppose that $
\kappa\ge\cf(\gamma)>\aleph_0$ and $A^f_i$ reflects on $
\gamma$. In particular, $A^f_i\cap e_\gamma$ is stationary in
$\gamma$. Given an $\alpha<\gamma$,
we can find a $\zeta\in A^f_i\cap e_\gamma$ such that $\alpha<\zeta$.
Then, there is a $\xi\in a^\zeta_i$ such that $\alpha<\xi$
and $f(\xi)\in a^\zeta_i$, by the definition of $A_i^f$.
Since $\cf(\gamma)\le\kappa\le\mu_i$, we have that
$a^\zeta_i\subseteq a^\gamma_i$, and we are done.

(C${}^{'}$) This follows immediately by (C).

(D) Define $h:S_{<\kappa}^\lambda\to\kappa$ by $h(\alpha)\deq\Min\{i:\,
\alpha\in A^f_i\}.$ Then, if $\delta\in S_f\cap
\acc(E)$, the function $h$ is defined on $S^\lambda_{<\kappa}
\cap\delta$ and not constant on any
stationary set of $\delta$, by (C${}^{'})$. By 1.11, $h$ is increasing on some
club of $\delta$.\eop${}_{2.7.}$

\medskip

{\bf Theorem 2.8.} (1) Suppose
that $\lambda=\mu^+$, and $\mu>\cf(\mu)=\kappa>\aleph_0$.

{\it Then,\/}

$\clubsuit^\ast_{-\mu}
(\lambda\setminus S^\lambda_\kappa)$ holds and
$\clubsuit^\ast_{-\mu}(S^\lambda_\kappa)/\II[\lambda,\kappa)$ holds.

(2) If in (1) we in addition assume that $\mu$ is a strong limit,
{\it then,\/}

$\club^\ast(\lambda\setminus S^\lambda_\kappa)$ holds and
$\club^\ast(S^\lambda_\kappa)/\II[\lambda,\kappa)$ holds.

{\bf Proof.} Let us fix sets $a_i^\alpha\,(i<\kappa)$ for $\alpha<\lambda$
as guaranteed by Lemma 2.7.
We fix for all $\alpha\in\lambda$ a cofinal subset $P_\alpha$ of $\alpha$
such that $\otp(P_\alpha)=\cf(\alpha).$

(1) We shall define for $\alpha\in\lambda$,
$$\PP_\alpha=\cases{ \{a^\alpha_i:\,i<\kappa
\,\&\,\sup(a^\alpha_i)=\alpha\}\cup\{P_\alpha\} & if $\card{\alpha}\ge\kappa$\cr
                     \{\alpha\}                 & otherwise.\cr}$$
Now, certainly each $\PP_\alpha$ is a family of $\le\card{\alpha}$
subsets of $\alpha$.

Suppose that $A\in[\lambda]^\lambda$ is given, and let $f$ be an increasing
enumeration of $A$. In particular, $f(\zeta)\ge\zeta$ for all $\zeta\in\lambda$.
The set
$$C\deq\{\alpha<\lambda:\,\forall\zeta\in\alpha\,(f(\zeta)<\alpha)\}\setminus
\kappa$$
is a club of $\lambda$.
Let us take any $\alpha\in C\cap S^\lambda_{\neq\kappa}$. By
(B) of Lemma 2.7, there is an $i<\kappa$ such that
$\alpha\in A^f_i$, where $A^f_i$ is as defined in Lemma
2.7.
Then
$\alpha=
\sup\{\zeta\in a^\alpha_i:\,f(\zeta)\in a^\alpha_i\}$, so
$\alpha=\sup(A\cap a^\alpha_i).$

This proves $\club^\ast_{-\mu}(\lambda\setminus S^\lambda_\kappa)$.
With the same definition of $\bar{\PP}\deq\langle\PP_\alpha:\,
\alpha<\lambda\rangle$, let us start again from
an $A\in[\lambda]^\lambda$, and $f$ and $C$ as above.
Let $S_f$ and $N$ be as in Lemma 2.7, so $S_f
\setminus N\in\II[\lambda,\kappa)$.
If $\alpha\in C\cap S^\lambda_\kappa\setminus S_f$,
then we argue as above, to conclude that $\alpha\in
G_{\clubsuit^\ast_{-\mu}(S^\lambda_\kappa)/\II[\lambda,\kappa)}^{\bar{\PP}}[A].$

(2) With the same notation as above, we define
$$\PP_\alpha\deq\cases{
\{\hbox{all cofinal sequences of
}\alpha\hbox{ included in } a^\alpha_i:\,i<\kappa\}
\cup\{P_\alpha\} & if $\card{\alpha}\ge\kappa$\cr
\{\alpha\} & otherwise.\cr}$$
Then $\card{\PP_\alpha}\leq 2^{\card{\alpha}}+\kappa<\mu$.

We argue similarly on $S^\lambda_\kappa$, using the set $S_f$ as above.\eop$
{}_{2.8.}$

\medskip

As a consequence, we obtain another proof of (a part of)
a theorem from [Sh 212, 14]
and [Sh 247, 6], as well as some other statements.

\medskip

{\bf Corollary 2.9.} Suppose that $\lambda=\mu^+$ and $\mu>
\cf(\mu)=\kappa>\aleph_0$. {\it Then\/}:

(1) No normal filter on $\lambda\setminus
 S^\lambda_\kappa$ is $\lambda^+$-saturated.

(2) If $\DD$ is a normal filter on $S^\lambda_\kappa$ such that
$\II[\lambda,\kappa)\rest S^\lambda_\kappa$ contains
only non-stationary sets, then $\DD$ is not $\lambda^+$-saturated.

(3) If $\club(S)$ holds, then $\DD_\lambda\rest S$ is not $\lambda^+$-saturated.

{\bf Proof.} (1) $\club^\ast_{-\mu}(\lambda\setminus S^\lambda_\kappa)$
holds, by 2.8. By 2.5, $\DD$ cannot be $\lambda^+$-saturated.

(2) Similar.

(3) See 2.10.2 below.\eop${}_{2.9.}$

\medskip

{\bf Concluding Remarks 2.10.} (0) Another useful version of the
$\club$-principle
on $\lambda$ for a stationary $S\subseteq\lambda$
is $\club^{-}(S)$, which says that

There is a sequence $\bar\PP=\langle \PP_\alpha:\,\alpha\in S\cap LIM \rangle$
such that $\PP_\alpha$ is a family of $\le\alpha$
unbounded subsets of $\alpha$
with the property that for all unbounded subsets $A$ of $\lambda$,
$$G_{\club^{-}(S)}^{\bar\PP}[A]\deq\{\alpha\in S:\,\exists B\in\PP_\alpha
(B\subseteq A\cap\alpha)\}$$
is non-empty.

As opposed to the situation with $\diamondsuit$, it
does not have to be true that
$\club^{-}(S)\implies\club(S)$. Of course, still $\club^\ast(S)\implies
\club^{-}(S)$, so we do not consider $\club^{-}$ here.

(1) We can also consider versions of $\club^\ast$
or $\club^{-}$ for which the size of 
each $\PP_\alpha$ is determined by some cardinal $\mu_\alpha$, not necessarily
equal to $\card{\alpha}$. Also, we can combine this idea with the idea
of $\club_{-\bar\theta}$, so also the order type of sets in $\PP_\alpha$ is
controlled by some prescribed sequence $\bar\theta$. 

(2) If we now define $\clubsuit^-_{-\mu}(S)$ in the
obvious way, then it follows from the proof of 2.5. that for $\lambda=
\mu^+$ and $\mu$ singular, $\club^-_{-\mu}(S)$ is enough to guarantee
that $\DD_\lambda\rest S$ is not $\lambda^+$-saturated. Therefore,
in particular, $\clubsuit(S)$ suffices.

For $\lambda$ the successor of a strong limit,  most ``reasonable''
versions of $\club$
coincide.

(3) After hearing
our lecture at the Logic Seminar in Jerusalem,
Fall 1994, M. Magidor showed us an alternative proof of 2.5 using
elementary embeddings and ultrapowers and not requiring $\mu$ to be
a limit cardinal.

(4) 
The assumptions of 1.13. and 2.5.
seem similar, but we
point out that there are in fact different.
The existence of a $\theta<\lambda$ which weakly reflects at $\kappa$
is not the same as the assumption that $\II[\lambda,\kappa)\rest
S^\lambda_\kappa$ contains only bounded sets.

\bigskip

\eject
\baselineskip=12pt

\centerline{\bf REFERENCES}
\bigskip

\item{[Ba]} J. E. Baumgartner,  A new class of order types,
{\it Annals of Mathematical Logic\/} (9) 187-222, 1976.

\item{[CDSh 571]} J. Cummings, M, D\v zamonja and S. Shelah, 
A consistency result on weak reflection, {\it Fundamenta
Mathematicae\/} (148) 91-100, 1995.

\item{[FMS]} M. Foreman, M. Magidor and S. Shelah,
Martin's maximum, saturated ideals and nonregular ultrafilters I,
{\it Annals of Mathematics, Second Series\/} (27) 1-47, 1988.

\item{[Fo]} M. Foreman, More saturated ideals, {\it Cabal Seminar 79-81,
Proceedings, Caltech-UCLA Logic Seminar 1979-81,
A.S. Kekris, D.A. Martin and Y.N. Moschovakis (eds.)\/},
Lecture Notes in Mathematics (1019) 1-27, {\it Springer-Verlag\/}

\item{[MkSh 367]} A.~H. Mekler and S. Shelah,
The consistency strength of ``every stationary set reflects'',
{\it Israel Journal of Mathematics,\/} (67) 353-366, 1989.

\item{[Ku 1]} K. Kunen, Saturated ideals, {\it Journal of
Symbolic Logic, \/} (43) 65-76, 1978.

\item{[Ku 2]} K. Kunen, Set theory, an Introduction to Independence
Proofs, {\it North-Holland,\/} Amsterdam 1980.

\item{[Os]} A .~J. Ostaszewski, On countably compact perfectly normal
spaces, {\it Journal of London Mathematical Society,\/} (14) 505-516, 1975.

\item{[Sh -g]} S. Shelah, Cardinal Arithmetic,
{\it Oxford University Press\/} 1994.

\item{[Sh 88a]} S. Shelah,
Appendix: on stationary sets ( to ``Classification of nonelementary
  classes. II. Abstract elementary classes''), in
{\it Classification theory (Chicago, IL, 1985),
Proceedings of the USA-Israel Conference on Classification Theory,
Chicago, December 1985; J.T. Baldwin (ed.) \/},
Lecture Notes in Mathematics,
(1292) 483--495, {\it Springer\/} Berlin 1987.

\item{[Sh 98]} S. Shelah, Whitehead
groups may not be free, even assuming $CH$,
{\it Israel Journal of Mathematics,\/} (35) 257-285, 1980.

\item{[Sh 108]} S. Shelah, On Successors of Singular Cardinals, In
{\it Logic Colloquium 78, M. Boffa, D. van Dalen, K. McAloon (eds.)\/}
357-380,
{\it North-Holland Publishing Company\/} 1979.

\item{[Sh 186]} S. Shelah,
Diamonds, Uniformization, {\it Journal of Symbolic
Logic \/} (49) 1022-1033, 1984.

\item{[Sh 212]} S. Shelah,
The existence of coding sets,
In {\it Around classification theory of models},
Lecture Notes in Mathematics, (1182) 188-202, {\it Springer}, Berlin
1986.

\item{[Sh 247]} S. Shelah,
More on stationary coding, In {\it Around classification theory of models},
Lecture Notes in Mathematics, (1182) 224-246, {\it Springer}, Berlin
1986.

\item{[Sh 351]} S.Shelah, Reflecting stationary sets and successors of
singular cardinals, {\it Archive for Mathematical Logic\/} (31) 25-53, 1991.

\item{[Sh 355]} S. Shelah, $\aleph_{\omega+1}$ has a Jonsson Algebra, in
{\it Cardinal Arithmetic,\/} Chapter II, {\it Oxford University
Press\/} 1994.

\item{[Sh 365]} S. Shelah, Jonsson Algebras in Inaccessible Cardinals, in
{\it Cardinal Arithmetic,\/} Chapter III, {\it Oxford University
Press\/} 1994.

\item{[Sh 420]} S. Shelah, Advances in Cardinal Arithmetic, In {\it
Finite and Infinite Combinatorics in Sets and Logic,} N. Sauer et al (eds.),
pp. 355-383, Kluwer Academic Publishers, 1993.

\item{[Sh 460]} S. Shelah, The Generalized Continuum Hypothesis revisited,
{\it Israel Journal of Mathematics, submitted.}

\item{[StvW]} J. R. Steel and R. van Wesep, Two Consequences of Determinacy
Consistent with Choice, {\it Transactions of the AMS,\/} (272)1, 67-85, 1982.

\baselineskip=24pt
\eject

\bye